\documentclass{article}
\RequirePackage{natbib}
\usepackage{amsmath, amssymb, amsthm, dsfont, xcolor, graphicx, times}
\usepackage[hidelinks]{hyperref}
\usepackage[margin=3.5cm]{geometry}

\newtheorem{theorem}{Theorem}
\newtheorem{corollary}{Corollary}
\newtheorem{lemma}{Lemma}
\newtheorem{prop}{Proposition}
\newtheorem{assumption}{Assumption}
\title{Only Closed Testing Procedures are Admissible for Controlling False Discovery Proportions}

\author{Jelle J. Goeman\footnote{Dept.\ of Biomedical Data Sciences, Leiden University Medical Center, Leiden, The Netherlands}
\and Jesse Hemerik\footnote{Oslo Centre for Biostatistics and Epidemiology, University of Oslo, Oslo, Norway}
\and Aldo Solari\footnote{Dept.\ of Economics, Management and Statistics, University of Milano-Bicocca, Milan, Italy}}

\date{Published in The Annals of Statistics: \url{http://doi.org/10.1214/20-AOS1999}}

\begin{document}

\maketitle

\begin{abstract}
We consider the class of all multiple testing methods controlling tail probabilities of the false discovery proportion, either for one random set or simultaneously for many such sets. This class encompasses methods controlling familywise error rate, generalized familywise error rate, false discovery exceedance, joint error rate, simultaneous control of all false discovery proportions, and others, as well as gene set testing in genomics and cluster inference in neuroimaging. We show that all such methods are either equivalent to a closed testing procedure, or are uniformly improved by one. Moreover, we show that a closed testing method is admissible if and only if all its local tests are admissible. This implies that, when designing  methods, it is sufficient to restrict attention to closed testing. We demonstrate the practical usefulness of this design principle by obtaining more informative inferences from the method of higher criticism, and by constructing a uniform improvement of a recently proposed method.
\end{abstract}

\section{Introduction}

Closed testing \citep{Marcus1976} is a fundamental principle of familywise error rate (FWER) control in multiple hypothesis testing. Indeed, almost every known procedure controlling FWER has been shown to be a special case of closed testing, and many procedures have been explicitly constructed as such. This is natural from a theoretical perspective, as  \cite{Sonnemann1982, Sonnemann2008}, and \cite{Sonnemann1988} have shown that closed testing is necessary for FWER control: every admissible procedure that controls FWER is a special case of closed testing. \cite{Romano2011} extended the results of Sonnemann and Finner, proving that from a FWER perspective not every closed testing procedure is admissible; only consonant procedures are. These results are valuable for designers of FWER controlling methods, who can rely exclusively on closed testing as a general design principle. Alternative design principles exist, such as the partitioning principle \citep{Finner2002} and sequential rejection \citep{Goeman2010}, but these are equivalent to closed testing.

Rather than only for FWER control, \cite{Goeman2011} showed that closed testing may also be used to obtain simultaneous confidence bounds for the false discovery proportion (FDP) of all subsets within a family of hypotheses. Used in this way, closed testing allows a form of post-selection inference. It allows users to look at the data prior to choosing thresholds and criteria for significance, while still keeping control of tail probabilities of the FDP. The approach of \cite{Goeman2011} is equivalent to an earlier approach by \cite{Genovese2004, Genovese2006} that did not explicitly use closed testing. A natural question that arises is whether similar results to those of \cite{Sonnemann1982}, \cite{Sonnemann1988} and \cite{Romano2011} also hold for this novel use of closed testing. When controlling FDP, is it sufficient to look only at closed testing-based methods? Which methods controlling FDP are admissible? These are the questions we will address in this paper.

\section{Overview and main results}

This paper has three main contributions. First, it presents a unification of methods and error rates, rewriting a wide range of diverse procedures as examples of a novel class. Within this class, we give first a necessary condition for admissibility, and then a  sufficient one. We start with a birds-eye view of these main results.

\cite{Genovese2004} and \cite{Goeman2011} considered simultaneous control of FDP for all subsets of a testing problem. For a family of hypotheses of interest $(H_i)_{i\in I}$, these authors have proposed methods to find upper $(1-\alpha)$-confidence bounds $\mathbf{q}^I(S)$ for the FDP $\pi_0(S)$, for all $S \subseteq I$, that are simultaneous for all such $S$. This means that
\begin{equation} \label{eq FDP control}
\mathrm{P}(\pi_0(S) \leq \mathbf{q}^I(S) \textrm{\ for all $S \subseteq I$}) \geq 1-\alpha.
\end{equation}
In this paper, we investigate the class of all methods controlling the error rate (\ref{eq FDP control}). In Section \ref{sec procedures} we show that many procedures that seem to target control of other quantities than FDP at first sight are members of this general class. These include all methods with regular FWER control; FWER control of intersection hypotheses; $k$-FWER control; simultaneous $k$-FWER control; False Discovery Exceedance control; control of the Joint Error Rate; and methods constructing confidence intervals for the overall proportion of true (or false) hypotheses. Essentially, the procedures we can rewrite as a special case of (\ref{eq FDP control}) are all procedures that control a tail probability of the number or proportion of false discoveries from above, either for one random set or simultaneously for several such sets.

Broad though this class of methods may be, it turns out that we can make strong statements that are valid for the whole class. We focus on admissibility, and therefore on the existence of uniform improvements of methods controlling (\ref{eq FDP control}). The first central result of this paper is given in Section \ref{sec thm suff}: given any method controlling (\ref{eq FDP control}), we can always construct a closed testing procedure that is either equivalent to the method we started with, or a uniform improvement of that method. Thus, our result implies that a necessary condition for admissibility is equivalence to a closed testing procedure. Moreover, we give an explicit construction of the improvement using the approach of \cite{Genovese2004} and \cite{Goeman2011}. The second main result of this paper is a sufficient condition for admissibility. We show in Section \ref{sec thm adm} that a closed testing procedure is admissible if the local tests that define the closed testing procedure are admissible.

Taken together, these results give design principles for multiple testing procedures. To design admissible procedures it is sufficient to create a closed testing procedure with admissible local tests. To show admissibility for a procedure designed in a different way, it is sufficient to show that the procedure is equivalent to such a procedure. We will discuss practical implications of our results for researchers seeking to develop new methods. We do this by revisiting two testing procedures: Higher Criticism \citep{Donoho2004, Meinshausen2006a} and the simultaneous FDP bounds of \cite{Katsevich2018}. In both cases we do not only uniformly improve the inferential statements of the methods, we also extend their scope by deriving non-trivial bounds for the FDP of sets these methods did not initially target.

\section{Inference on false discovery proportions} \label{sec procedures}

Assume that we have data $\mathbf{X}$ distributed according to some unknown probability distribution $\mathrm{P} \in \Omega$. About $\mathrm{P}$ we may formulate hypotheses of the form $H \subseteq \Omega$. Let the family of hypotheses of interest be $(H_i)_{i\in I}$, where $I \subseteq C \subseteq \mathbb{N}$ is finite. The set $C$, possibly infinite, is arbitrary here, but will become important in Section \ref{sec monotone}. Within the family $I$, let $I_0 = \{i \in I\colon \mathrm{P} \in H_i\}$ be the index set of the true hypotheses and $I_1 = I \setminus I_0$ the index set of the false hypotheses. We will make no further model assumptions in this paper: any models, any test statistics, and any dependence structures will be allowed. Equalities and inequalities between random variables should be read as holding almost surely for all $\mathrm{P} \in \Omega$ unless otherwise stated. Proofs of all theorems, lemmas and propositions are in Section E in the Supplemental Information.
Throughout the paper we will denote all random quantities in boldface. Upper case variables (except $\mathrm{P}$) always refer to sets.

We will be studying procedures with FDP control. The FDP of a finite set $S$ is given by
\[
\pi_0(S) = \frac{|S \cap I_0|}{|S| \vee 1}.
\]
We define a procedure with FDP control on $I$ (i.e.\ on $(H_i)_{i\in I}$) as a random function $\mathbf{q}^I\colon 2^I \to [0,1]$, where $2^I$ is the power set of $I$, such that for all $\mathrm{P} \in \Omega$ it satisfies (\ref{eq FDP control}).

It will be more convenient to use an equivalent representation that gives a simultaneous lower $(1-\alpha)$-confidence bound for $|S \cap I_1|$, the number of true discoveries. We say that a random function $\mathbf{d}^I\colon 2^I \to \mathbb{R}$ has a $(1-\alpha)$-\emph{true discovery guarantee} on $I$ if, for all $\mathrm{P} \in \Omega$,
\begin{equation} \label{eq discovery guarantee}
\mathrm{P}(\mathbf{d}^I(S) \leq |S \cap I_1|  \textrm{\ for all $S \subseteq I$}) \geq 1-\alpha.
\end{equation}
We will usually suppress the dependence on $\alpha$ when talking about true discovery guarantees. To see that the class of methods of FDP control and the class of procedures with a true discovery guarantee are equivalent, note that if $\mathbf{q}^I$ fulfils (\ref{eq FDP control}), then
\[
\mathbf{d}^I(S) = (1-\mathbf{q}^I(S))|S|,
\]
fulfils (\ref{eq discovery guarantee}) and, if $\mathbf{d}^I$ fulfils (\ref{eq discovery guarantee}), then
\[
\mathbf{q}^I(S) = \frac{|S|-\mathbf{d}^I(S)}{|S| \vee 1}
\]
fulfills (\ref{eq FDP control}). In the rest of the paper we will focus on true discovery guarantee procedures, which are mathematically easier to work with than methods with FDP control, e.g.\ because they automatically avoid issues with empty sets $S$. Without loss of generality we may assume that $\mathbf{d}^I(S)$ takes integer values, and that $0 \leq \mathbf{d}^I(S) \leq |S|$. If $\mathbf{d}^I(S)$ is not integer, we may freely replace $\mathbf{d}^I(S)$ by $\lceil \mathbf{d}^I(S) \rceil$.

The class of FDP control (cf.\ true discovery guarantee) procedures encompasses seemingly diverse methods. Only few authors \citep{Genovese2006, Goeman2011,  Goeman2017, Blanchard2017} have explicitly proposed procedures that target control of FDP for all sets $S$ simultaneously as implied by (\ref{eq FDP control}). However, many other well-known types of multiple testing procedures turn out to be special cases of FDP control procedures, even if they were not directly formulated to control (\ref{eq FDP control}) or its equivalent. We will review these procedures briefly in the rest of this section in order to emphasize the wide range of applications of the results of this paper. We will reformulate such procedures in terms of $\mathbf{d}^I$.

Procedures that control FWER \citep[e.g.][]{Westfall1993, Bretz2009, Berk2013, Janson2016} within the family defined by $I$ are usually defined as producing a random set $\mathbf{K}$ (possibly empty) for which it is guaranteed that, for all $\mathrm{P} \in \Omega$,
\(
\mathrm{P}(|\mathbf{K} \cap I_0| = 0) \geq 1-\alpha,
\)
A generalization, $k$-FWER \citep{Hommel1988, Lehmann2005, Romano2006, Sarkar2007,  Guo2010, Finos2011}, makes sure that, for all $\mathrm{P} \in \Omega$,
\[
\mathrm{P}(|\mathbf{K}\cap I_0| < k) \geq 1-\alpha,
\]
which reduces to regular FWER if $k=1$ is chosen.
It is easily seen that this is equivalent to requiring (\ref{eq discovery guarantee}) if we take
\begin{equation} \label{eq kFWER}
\mathbf{d}^I(S) = \left\{ \begin{array}{ll} |S|-k+1 & \textrm{if $S = \mathbf{K}$,} \\ 0 & \textrm{otherwise.} \end{array} \right.
\end{equation}
Free additional statements may be obtained from (\ref{eq kFWER}) by direct logical implication. For example, if $\mathbf{d}^I(S) = |S|-k+1$ then we may immediately set $\mathbf{d}^I(U) = |U|-k+1$, if positive, for all $U \subseteq S$ without compromising (\ref{eq discovery guarantee}). We will come back to such implications in Section \ref{sec coherence}.

Related to $k$-FWER are methods controlling False Discovery Exceedance (FDX), also known as $\gamma$-FDP, at level $\boldsymbol\gamma$ \citep{Dudoit2004, Korn2004, Romano2006, Farcomeni2009, Sun2015, Delattre2015}. Such methods find a random set $\mathbf{K}$ (possibly empty) such that, for all $\mathrm{P} \in \Omega$,
\[
\mathrm{P}\Big(|\mathbf{K} \cap I_0| \leq \boldsymbol\gamma|\mathbf{K}|\Big) \geq 1-\alpha,
\]
which is equivalent to (\ref{eq discovery guarantee}) with
\[
\mathbf{d}^I(S) = \left\{ \begin{array}{ll} \big\lceil(1-\boldsymbol\gamma)|S|\big\rceil & \textrm{if $S = \mathbf{K}$,} \\ 0 & \textrm{otherwise.} \end{array} \right.
\]
In most methods controlling FDX the control level $\boldsymbol\gamma$ is fixed, but it may also be random, as e.g.\ in the permutation-based method of \cite{Hemerik2018}. Variants, such as kFDP \citep{Guo2014}, which allow a minimum number of false discoveries regardless of the size of $\mathbf{K}$, also fit (\ref{eq discovery guarantee}).

Other methods allow $\gamma$ to be chosen post-hoc by controlling FDX simultaneously over several values of $\gamma$. One way to achieve this is by control of the Joint Error Rate (JER). The JER \citep{Blanchard2017} constructs a sequence of $\mathbf{m}\geq 0$ distinct random sets $\mathbf{K}_1,\ldots,\mathbf{K}_\mathbf{m} \subseteq I$ and corresponding random bounds $\mathbf{k}_1, \ldots, \mathbf{k}_\mathbf{m}$, such that, for all $\mathrm{P} \in \Omega$,
\[ 
\mathrm{P}\big(|\mathbf{K}_i\cap I_0| < \mathbf{k}_i \textrm{\ for all $1\leq i \leq \mathbf{m}$}\big) \geq 1-\alpha.
\]
This is a special case of (\ref{eq discovery guarantee}) if we set
\[
\mathbf{d}^I(S) = \left\{ \begin{array}{ll} |\mathbf{K}_i| - \mathbf{k}_i+1 & \textrm{if $S = \mathbf{K}_i$ for some $1\leq i \leq \mathbf{m}$,} \\ 0 & \textrm{otherwise.} \end{array} \right.
\]
Joint error rate control may be used with nested sets \citep{Blanchard2017} or tree-structured sets \citep{Durand2018}, and is meant to be combined with interpolation (see Section \ref{sec coherence}). Similar approaches were used by e.g.\ the permutation-based methods of \cite{Meinshausen2006} and \cite{Hemerik2018a}. Also the approach of \cite{Katsevich2018}, discussed in detail in Section \ref{sec KR}, can be seen as controlling JER with nested sets.

A different category of methods involves FWER control of many intersection hypotheses, as e.g.\ used in gene set testing in genomics and in cluster inference in neuroimaging. In genomics, a collection of distinct sets $K_1, \ldots, K_m \subseteq I$ is given a priori, and the procedure generates corresponding random indicators $\mathbf{k}_1, \ldots, \mathbf{k}_m \in \{0,1\}$ for detection of signal in the corresponding set. FWER is controlled over all statements made, i.e., for all $\mathrm{P} \in \Omega$,
\begin{equation} \label{eq DAG FWER}
\mathrm{P}\big(|K_i\cap I_1| \geq \mathbf{k}_i \textrm{\ for all $i=1,\ldots, m$}\big) \geq 1-\alpha.
\end{equation}
This corresponds to (\ref{eq discovery guarantee}) with
\[
\mathbf{d}^I(S) = \left\{ \begin{array}{ll} \mathbf{k}_i & \textrm{if $S = K_i$ for some $1\leq i \leq m$,} \\ 0 & \textrm{otherwise.} \end{array} \right.
\]
Examples of such methods include \cite{Meinshausen2008}, \cite{Goeman2008}, \cite{Goeman2012}, \cite{Meijer2015a}, \cite{Meijer2015b}, and \cite{Meijer2015}. In the latter two papers a connection with FDP control was already noted. In neuroimaging, cluster inference methods are similar except that in this case the sets $\mathbf{K}_1, \ldots, \mathbf{K}_\mathbf{m}$ and their number $\mathbf{m}\geq 0$ are random, and $k_i=1$ for $i=1,\ldots,\mathbf{m}$ is fixed \citep{Poline1993}. FWER control (\ref{eq DAG FWER}) is guaranteed by Gaussian random field theory. Such control translates to a true discovery guarantee (\ref{eq discovery guarantee}) in the same way.

In partial conjunction testing \citep{Benjamini2008, Wang2018}, the hypothesis $H_0^{k/n}: |I_1| < k$ is tested for some $1\leq k \leq n$. The requirement that $\boldsymbol\delta$, taking values in $\{0,1\}$ is a valid test of $H_0^{k/n}$ is equivalent to (\ref{eq discovery guarantee}) with
\[
\mathbf{d}^I(S) = \left\{ \begin{array}{ll} \boldsymbol\delta k & \textrm{if $S = I$,} \\ 0 & \textrm{otherwise.} \end{array} \right.
\]
Finally, related to partial conjunction methods are methods that aim to make one-sided confidence intervals for $\pi_0(I)$, the proportion of true null hypotheses in the testing problem as a whole \citep{Meinshausen2006a, Ge2012}. Here, the requirement that $[0, \mathbf{u}]$ is a valid confidence interval for $\pi_0(I)$ is equivalent to demanding (\ref{eq discovery guarantee}) with
\[
\mathbf{d}^I(S) = \left\{ \begin{array}{ll} (1-\mathbf{u})|I| & \textrm{if $S = I$,} \\ 0 & \textrm{otherwise.} \end{array} \right.
\]

This listing of the different types of methods that may be written as true discovery guarantee methods is certainly not exhaustive, but a general pattern emerges. Any method controlling a $(1-\alpha)$-tail probability of the number or proportion of true discoveries (from below) or false discoveries (from above) either in one subset of $I$, or in several subsets simultaneously, are special cases of general discovery control procedures. The sets and bounds are all allowed to be random; only $\alpha$ must be fixed.

Writing procedures as true discovery guarantee procedures, even when the rewriting is trivial, may bring a new perspective to the use of the procedure. As proposed by \cite{Goeman2011}, procedures that fulfil (\ref{eq FDP control}) or (\ref{eq discovery guarantee}) allow a different, flexible way of using multiple testing methods. In flexible multiple testing the user may look at the data before choosing post hoc one or several sets $S \subseteq I$ of interest, based on any desired criteria, and find their $\mathbf{d}^I(S)$. Regardless of this data peeking the bounds on the selected sets are simultaneously valid due to the simultaneity in (\ref{eq discovery guarantee}). Writing procedures in this form, therefore, in principle opens the way to their use as post-selection inference methods \citep[see][for applications]{Rosenblatt2018, Ebrahimpoor2019}. Of course, this is only useful if the user has some real choice, i.e.\ if  $\mathbf{d}^I(S) \neq 0$ for a number of sets $S$. We will see in Section \ref{sec coherence} how to get rid of some of the zeros in the definitions above.

\section{True discovery guarantee using closed testing} \label{sec CT}

A general way to construct true discovery guarantee procedures is provided by closed testing, introduced by \cite{Marcus1976} for FWER control. \cite{Genovese2006} and \cite{Goeman2011} adapted closed testing to make it usable for true discovery guarantee and FDP control. We will briefly review these methods here.

For every finite set $S \subseteq C$ we define a corresponding intersection hypothesis as $H_S = \bigcap_{i \in S} H_i$. This hypothesis is true if and only if all $H_i$, $i \in S$ are true. We have $H_\emptyset = \Omega$, which is always true. For every intersection hypothesis $H_S$ we may choose a \emph{local test} $\boldsymbol\phi_S$, taking values in $\{0,1\}$, with 1 indicating rejection of $H_S$. This is a valid statistical test for $H_S$ if it has the property that, for all $\mathrm{P} \in H_S$
\[
\mathrm{P}(\boldsymbol\phi_S = 1) \leq \alpha.
\]
We always choose $\boldsymbol\phi_\emptyset = 0$ surely. Choosing a local test for every finite $S \subseteq C$ will yield a \emph{suite of local tests} $\boldsymbol\phi = (\boldsymbol\phi_S)_{S \subseteq C, |S|<\infty}$. To deal with restricted combinations \citep{Shaffer1986} efficiently, if present, we demand that identical hypotheses have identical tests: if for $U, V \subseteq C$ we have $H_U=H_V$, then $\boldsymbol\phi_U=\boldsymbol\phi_V$. If $H_U = \emptyset$ for some $U \subseteq C$, we may take $\boldsymbol\phi_U=1$ surely.

From a suite of local tests we may obtain a true discovery guarantee procedure in two simple steps. First, we need to correct the tests for multiple testing. We define the \emph{effective local test} within the family $I$ by
\[
\boldsymbol{\phi}_S^I = \min\{\boldsymbol\phi_U\colon S \subseteq U \subseteq I\}.
\]
As shown by \cite{Marcus1976}, the effective local tests have FWER control over all intersection hypotheses $H_S$, $S \subseteq I$, i.e., for all $\mathrm{P} \in \Omega$,
\[
\mathrm{P}(\boldsymbol{\phi}_S^I \leq |S \cap I_1| \textrm{\ for all $S \subseteq I$}) \geq 1-\alpha.
\]
Next, we calculate $\mathbf{d}^I(S)$. We see that the procedure defined by $\mathbf{d}^I(S) = \boldsymbol\phi_S^I$ already fulfils (\ref{eq discovery guarantee}). More recently, however, \cite{Goeman2011} showed that closed testing may also be used for more powerful FDP control. For any suite of local tests $\boldsymbol\phi$, these authors defined the associated procedure
\begin{equation} \label{def d}
\mathbf{d}^I_{\boldsymbol\phi}(S) = \min_{U \in 2^S}\{|S\setminus U|\colon \boldsymbol{\phi}^I_U = 0\},
\end{equation}
and proved the true discovery guarantee. Note that the minimum is always defined since $\boldsymbol{\phi}^I_\emptyset = \boldsymbol\phi_\emptyset = 0$ surely.

An earlier general approach to developing true discovery guarantee procedures was developed, without reference to closed testing, by \cite{Genovese2004, Genovese2006}. Starting from a suite of local tests, they proved coverage for the general true discovery guarantee procedure

\begin{equation} \label{def g}
\mathbf{g}^I_{\boldsymbol\phi}(S) = \min_{V \in 2^I} \{|S \setminus V|\colon \boldsymbol\phi_V =0\}.
\end{equation}

The difference between approaches (\ref{def d}) and (\ref{def g}) is that (\ref{def d}) uses a two-step approach, first correcting the local tests for multiple testing using the closed testing procedure, while (\ref{def g}) works directly on the local tests. In compensation, (\ref{def d}) only needs to look through the subsets of the set of interest $S$, while (\ref{def g}) looks through all subsets of the family $I$. The end result, however, is identical \citep{Hemerik2018a}:

\begin{lemma} \label{thm GW}
$\mathbf{g}^I_{\phi} = \mathbf{d}^I_{\phi}$.
\end{lemma}

The expressions (\ref{def d}) and (\ref{def g}) are very useful for constructing true discovery guarantee procedures. Local tests tend to be easy to specify in most models, as each local test is a test of a single hypothesis, so that standard statistical test theory may be used. Given a suite of local tests, (\ref{def d}) or (\ref{def g}) takes care of the multiplicity. A computational problem remains: direct application of (\ref{def d}) or (\ref{def g}) takes exponential time. Often, however, shortcuts are available that allow faster computation \citep{Goeman2011, Goeman2017, Dobriban2018}. We'll see examples in Sections \ref{sec HC} and \ref{sec KR}.

Comparing (\ref{def d}) and (\ref{def g}), the single step expression of \cite{Genovese2006} is clearly more elegant. However, the link of (\ref{def d}) to closed testing is valuable because it connects true discovery guarantee procedures to the enormous literature on closed testing \citep[see][for an overview]{Henning2015}. The detour via effective local tests is often profitable in practice because expressions for $\mathbf{d}^I_{\boldsymbol\phi}(S)$ can be easier to derive through expressions for $\boldsymbol\phi_S^I$ \citep{Hemerik2018, Goeman2017}.

\section{Coherence and interpolation} \label{sec coherence}

By viewing methods in terms of true discovery guarantees, as we have done in Section \ref{sec procedures}, they are upgraded from making a confidence statement about discoveries in a limited number of sets $S \subseteq I$ to doing the same for all subsets of $I$. However, in the definitions of Section \ref{sec procedures}, most of these statements are the trivial $\mathbf{d}^I(S)=0$. Often, however, some of the statements can be uniformly improved by a process called \emph{interpolation}. In this section we discuss interpolation and how it can improve true discovery guarantee procedures. We will define coherent procedures as procedures that cannot be improved by interpolation.

Let $\mathbf{d}^I$ be some true discovery guarantee procedure. We define the \emph{interpolation} $\mathbf{\bar d}^I$ of $\mathbf{d}^I$ as
\begin{equation} \label{def interpolation}
\mathbf{\bar d}^I(S) = \max_{U \in 2^I} \Big\{ \mathbf{d}^I(U) - |U \setminus S| + \mathbf{d}^I(S \setminus U) \Big\}.
\end{equation}
Interpolation was used in weaker versions or in specific cases by several authors \citep{Genovese2006, Meinshausen2006, Blanchard2017, Durand2018}. Taking $U = S$, we see that $\mathbf{\bar d}^I(S) \geq \mathbf{d}^I(S)$. Moreover, the improvement from $\mathbf{d}^I$ to $\mathbf{\bar d}^I$ is for free, as noted in the following lemma.

\begin{lemma} \label{thm interpolation}
If $\mathbf{d}^I$ is a true discovery guarantee procedure then so is $\mathbf{\bar d}^I$.
\end{lemma}

Intuitively, the rationale for interpolation is as follows. If $\mathbf{d}^I(U)$ is large, and $S$ has so much overlap with $U$ that the signal $\mathbf{d}^I(U)$ in $U$ does not fit in $U \setminus S$, then the remaining signal must be in $S$. Since this reasoning follows by direct logical implication, it will not increase the occurrence of type I error: we can only make an erroneous statement about $S$ if we had already made one about $U$.
As an example, consider interpolation for $k$-FWER controlling procedures. The interpolated version of (\ref{eq kFWER}) is simply
\begin{equation} \label{eq interpolate FWER}
\mathbf{\bar d}^I(S) =  0 \vee (|S \cap \mathbf{K}|-k+1),
\end{equation}
an expression that simplifies even further to $\mathbf{\bar d}^I(S) = |S \cap \mathbf{K}|$ with regular FWER when $k=1$.

Interpolation is not necessarily a one-off process, and interpolated procedures may sometimes be further improved by another round of interpolation. We call a procedure \emph{coherent} if it cannot be improved by interpolation, i.e.\ if
\begin{equation} \label{def_coh}
\mathbf{\bar d}^I(S) = \mathbf{d}^I(S) \textrm{\quad for all $S \subseteq I$}.
\end{equation}
We can characterize coherent procedures further with the following lemma.

\begin{lemma} \label{thm coherence}
$\mathbf{d}^I$ is coherent if and only if for every disjoint $V,W \subseteq I$ we have
\[ 
\mathbf{d}^I(V) + \mathbf{d}^I(W) \leq \mathbf{d}^I(V \cup W) \leq \mathbf{d}^I(V) + |W|.
\]
\end{lemma}

We intentionally use the same term \emph{coherent} that was used by \cite{Sonnemann1982} in the context of FWER control of intersection hypotheses. Looking only at FWER control of intersection hypotheses is equivalent to looking only at $\mathds{1}\{\mathbf{d}^I(S) >0\}$ for every $S$, where $\mathds{1}\{\cdot\}$ denotes an indicator function. In that case (\ref{def_coh}) reduces to simply requiring that $U \subseteq V$ and $\mathbf{d}^I(U) >0$ implies that $\mathbf{d}^I(V) >0$, which is exactly Sonnemann's definition of coherence.

Methods that are created through closed testing are automatically coherent, as the following lemma claims.

\begin{lemma} \label{thm CT coherent}
The procedure $\mathbf{d}_\phi^I$ is coherent.
\end{lemma}

Since an incoherent procedure can always be replaced by a coherent procedure that is at least as good, we will restrict attention to coherent procedures for the rest of this paper.

\section{Monotone procedures} \label{sec monotone}

The methods from the literature discussed in Sections \ref{sec procedures} and \ref{sec CT} are usually not defined for a specific family $I$ of hypotheses, but as generic procedures that can be used for any family, large or small. Researchers developing methods are usually not looking for good properties for a specific family at a specific scale $|I|$, but for methods that are generally applicable and have good properties whatever $I$.

We can embed the procedure $\mathbf{d}^I$ into a stack of procedures $\mathbf{d} = (\mathbf{d}^I)_{I \subseteq C, |I|<\infty}$, where we may have some maximal family $C \subseteq \mathbb{N}$. We will briefly call $\mathbf{d}$ a \emph{monotone procedure} if it fulfils the three criteria below. In contrast, we call $\mathbf{d}^I$ for a specific $I$ a \emph{local procedure}, or a local member of $\mathbf{d}$.
\begin{enumerate}
\item \emph{true discovery guarantee}: $\mathbf{d}^I$ is a true discovery guarantee procedure for every finite $I \subseteq C$;
\item \emph{coherence}: $\mathbf{d}^I$ is coherent for every finite $I \subseteq C$;
\item \emph{monotonicity}: $\mathbf{d}^I(S) \geq \mathbf{d}^J(S)$ for every finite $S \subseteq I \subseteq J \subseteq C$.
\end{enumerate}
The first two criteria are no more than natural. We demand a true discovery guarantee for every member of the monotone procedure, and we demand coherence for every local member since otherwise we may always improve it by a coherent procedure. The monotonicity requirement relates local procedures at different scales to each other. It says that inference on the number of discoveries in a set $S$ should never get better if we embed $S$ in a larger family $J$ rather than in a smaller family $I$. As the multiple testing problem gets larger, inference should get more difficult. This requirement relates closely to the ``subsetting property'' of \cite{Goeman2014} and the monotonicity property of various FWER control procedures \citep[e.g.][]{Bretz2009, Goeman2010}. It is a natural requirement, and the procedures cited in Section \ref{sec procedures} generally adhere to it by construction.

There are a few notable exceptions to the rule that method designers tend to design monotone rather than local procedures. All the examples we are aware of are FWER-controlling procedures. \cite{Rosenblum2014} proposed a local procedure for $|I|=2$ hypotheses that optimizes the power for rejecting at least one of these. Their method is specific for the scale $|I|$ it was defined for; extensions to $|I|>2$ do not exist \citep{Rosset2018}. In another example, \cite{Rosset2018} developed methods that optimize the power for detecting at least one true effect for specific scales $|I|$ under an exchangeability assumption. These methods also have non-monotone behavior.

We remark, however, that every coherent local $\mathbf{d}^I$ true discovery guarantee procedure may be trivially embedded in a monotone procedure with $C=I$ (or even $C=\mathbb{N}$) by setting
\begin{equation} \label{eq trivial scalability}
\mathbf{d}^J(S) = \left\{ \begin{array}{ll} \mathbf{d}^I(S) & \textrm{if $S \subseteq I$,} \\ 0 & \textrm{otherwise.} \end{array} \right.
\end{equation}
This embedding allows translation of properties of monotone procedures to properties of their local members. We will mostly be studying monotone procedures in this paper, but investigate implications for local procedures where appropriate.

Procedures created using closed testing are automatically monotone, as formalized in the following lemma.

\begin{lemma} \label{thm CT monotone}
The procedure $\mathbf{d}_{\boldsymbol\phi} = (\mathbf{d}_{\boldsymbol\phi}^I)_{I \subseteq C, |I|<\infty}$ is a monotone procedure.
\end{lemma}

The property of primary interest to us is admissibility. Let us formally define admissibility for true discovery guarantee procedures. Recall that a statistical test $\boldsymbol\delta$ of a hypothesis $H$ is uniformly improved by a statistical test $\boldsymbol{\tilde\delta}$ of the same hypothesis if (1.) $\boldsymbol{\tilde\delta} \geq \boldsymbol{\delta}$; and (2.) $\mathrm{P}(\boldsymbol{\tilde\delta} > \boldsymbol\delta)>0$ for some $\mathrm{P} \in \Omega$. A statistical test is \emph{admissible} if no test exists that uniformly improves it \citep[Section 6.7]{Lehmann2006}. We call a suite of local tests $\boldsymbol\phi$ admissible if $\boldsymbol\phi_S$ is admissible for all finite $\emptyset \subset S \subseteq C$. We note that existence of admissible tests is not assured in all models, but that under a weak condition all tests that exhaust the $\alpha$-level are admissible. We discuss these technical issues in Section A in the Supplemental Information, where we also motivate our definition of admissibility compared to alternatives in the literature.

Analogously to admissibility of single tests we define admissibility for true discovery guarantee procedures. A uniform improvement of a monotone procedure $\mathbf{d}$ is a monotone procedure $\mathbf{\tilde{d}}$ such that (1.) $\mathbf{\tilde{d}}^I(S) \geq \mathbf{d}^I(S)$ for all finite $S \subseteq I \subseteq C$; and (2.) $\mathrm{P}(\mathbf{\tilde{d}}^I(S) > \mathbf{d}^I(S)) > 0$ for some $\mathrm{P} \in \Omega$ and some finite $S \subseteq I \subseteq C$. A uniform improvement of a local procedure $\mathbf{d}^I$ is a local procedure $\mathbf{\tilde{d}}^I$ such that (1.) $\mathbf{\tilde{d}}^I(S) \geq \mathbf{d}^I(S)$ for all $S \subseteq I$; and (2.) $\mathrm{P}(\mathbf{\tilde{d}}^I(S) > \mathbf{d}^I(S)) > 0$ for some $\mathrm{P} \in \Omega$ and some $S \subseteq I$. We call a local or monotone procedure that cannot be uniformly improved admissible. If all local members of a monotone procedure are admissible, then the monotone procedure is admissible, but the converse is not necessarily true, as illustrated in Section B in the Supplemental Information.

\section{All admissible procedures are closed testing procedures} \label{sec thm suff}

Theorem \ref{thm sufficiency}, below, claims that every monotone true discovery guarantee procedure is either equivalent to a closed testing procedure or can be uniformly improved by one. We already know from Lemma \ref{thm coherence} that every incoherent procedure can be uniformly improved by a coherent procedure. It follows that every procedure that is not equivalent to a closed testing procedure is inadmissible: the class of all closed testing procedures is essentially complete \citep[Section 1.8]{Lehmann2006} for procedures with a true discovery guarantee, and therefore for FDP control. This is the first main result of this paper.

\begin{theorem} \label{thm sufficiency}
Let $\mathbf{d}$ be a monotone procedure. Then, for every finite $S\subseteq C$,
\[
\boldsymbol\phi_S = \mathds{1}\{\mathbf{d}^{S}(S) > 0\}
\]
is a valid local test of $H_S$. For the suite $\boldsymbol\phi = (\boldsymbol\phi_S)_{S\subseteq C, |S|<\infty}$ we have, for all $S \subseteq I\subseteq C$ with $|I|<\infty$,
\[
\mathbf{d}^I_{\boldsymbol\phi}(S) \geq \mathbf{d}^I(S).
\]
\end{theorem}

Coherence is necessary but not sufficient to guarantee admissibility. The procedure $\mathbf{d}^I_{\boldsymbol\phi}(S)$ implied by Theorem \ref{thm sufficiency} may in some cases be truly a uniform improvement over the original, coherent $\mathbf{d}^I(S)$.  To see a classical example in which a coherent procedure can uniformly improved by closed testing, think of Bonferroni. Combined with (\ref{eq interpolate FWER}), Bonferroni is coherent. However, it is uniformly improved by Holm's procedure that follows from a well-known step-down argument that incorporates an estimate of $\pi_0(I)$ into the procedure. This stepping-down can be seen as a direct application of closed testing with the local test defined in Theorem \ref{thm sufficiency}. Step-down arguments are standard for FWER control and have been applied to several FDP controlling methods in the past \citep{Blanchard2017, Goeman2017, Hemerik2018a}.

It should be noted that in case of a monotone procedure, the local test $\boldsymbol\phi_S$ defined in Theorem \ref{thm sufficiency} is truly local, in the sense that it uses only the information used by the restricted testing problem $\mathbf{d}^S$ about the hypotheses $H_i$, $i \in S$. For example, in a testing problem based on $p$-values, the local test would use only the $p$-values $p_i$, $i \in S$. In other testing problems, some global information may be used, e.g. the overall estimate of $\sigma^2$ in a large one-way ANOVA, but still in such situations the local test is very natural: as a local test for $H_S$ we use the test for discovery of signal in hypotheses $H_i$, $i\in S$, that we would use in the situation where the hypotheses $H_i$, $i\notin S$ are not of interest to us. Such a local test is implicitly defined by the local procedure $\mathbf{d}^S$.

The result of the theorem is formulated in terms of monotone procedures. It applies immediately to local procedures as well if we use the trivial embedding (\ref{eq trivial scalability}) of a local procedure into a monotone one. With this embedding we even have $\mathbf{d}^I_{\boldsymbol\phi}(S) = \mathbf{d}^I(S)$. This leads to the following corollary.

\begin{corollary} \label{thm corol}
Let $\mathbf{d}^I$ be a coherent procedure. Then, for every $S\subseteq I$,
\[
\boldsymbol\phi_S = \mathds{1}\{\mathbf{d}^{I}(S) > 0\}
\]
is a valid local test of $H_S$. For the suite $\boldsymbol\phi = (\boldsymbol\phi_S)_{S\subseteq I}$ we have, for all $S \subseteq I$,
\[
\mathbf{d}^I_{\boldsymbol\phi}(S) = \mathbf{d}^I(S).
\]
\end{corollary}

Corollary \ref{thm corol} shows that every coherent true discovery guarantee procedure is equivalent to a closed testing procedure. It may possibly be uniformly improved by another closed testing procedure if the suite of local tests $\boldsymbol\phi$ is not admissible, as we shall see in the next section.

Corollary \ref{thm corol} also confirms the equivalence between the closed testing and partitioning principles for FWER control. This has been clear since \cite{Finner2002} showed that closed testing procedures may be rewritten as partitioning procedures and that this sometimes uniformly improves them, while \cite{Sonnemann1982} and \cite{Sonnemann1988} had already shown that the family of closed testing procedure is complete for FWER control. However, since the result is important and, as far as we know, not explicitly stated in the literature we phrase it as a separate theorem.

\begin{theorem} \label{thm equiv_prin}
For every closed testing procedure there exists a partitioning procedure that rejects exactly the same hypotheses. For every partitioning procedure there exists a closed testing procedure that rejects exactly the same hypotheses.
\end{theorem}

Since both closed testing and partitioning procedures may be written as sequential rejection procedures \citep{Goeman2010}, while it cannot improve upon them by Corollary \ref{thm corol} and Theorem \ref{thm equiv_prin}, sequential rejection could be labelled a third equivalent principle.

\section{All closed testing procedures are admissible} \label{sec thm adm}

So far we have seen that a true discovery guarantee procedure may be uniformly improved by interpolation to coherent procedures, which in turn may be uniformly improved by closed testing procedures. Clearly, equivalence to a closed testing procedure is necessary for admissibility. Are all closed testing procedures admissible? In this section we derive a simple condition for admissibility of monotone procedures that is both necessary and sufficient. We show that admissibility of the monotone procedure $\mathbf{d}_{\boldsymbol\phi}$ follows directly from admissibility of its local tests. This is the second main result of this paper.

\begin{theorem} \label{thm admissibility}
$\mathbf{d}_{\boldsymbol\phi}$ is admissible if and only if the suite $\boldsymbol\phi$ is admissible.
\end{theorem}

We have already seen from Theorem \ref{thm sufficiency} that only closed testing procedures are admissible. Theorem \ref{thm admissibility} says that all closed testing procedures are admissible, provided they fulfil the reasonable demand that they are built from admissible local tests. To check admissibility of the local tests, Section A in the Supplemental Information shows that under a weak assumption it is sufficient to check that the local tests exhaust the $\alpha$-level. Theorem \ref{thm admissibility} thus makes it easy to guarantee admissibility of monotone procedures.

Unlike Theorem \ref{thm sufficiency}, the result of Theorem \ref{thm admissibility} does not immediately translate to local procedures: even if $\boldsymbol\phi$ is admissible, it may happen for some finite $I \subseteq C$ that $\mathbf{d}^I_{\boldsymbol\phi}$ can be uniformly improved by some other procedure $\mathbf{d}^I$. About such local improvements we have the following proposition.

\begin{prop} \label{thm locally admissible}
If $\mathbf{d}^I \geq \mathbf{d}^I_{\boldsymbol\phi}$ is admissible, then there is an admissible $\boldsymbol\psi$ such that $\mathbf{d}^I = \mathbf{d}^I_{\boldsymbol\psi}$ and, for all $S \subseteq I$, $\boldsymbol\psi_S \geq \boldsymbol\phi^I_{S}$.
\end{prop}

Proposition \ref{thm locally admissible} limits the available room for local improvements of admissible monotone procedures. Combining Proposition \ref{thm locally admissible} and Theorem \ref{thm admissibility} we see that such improvements have to be admissible monotone procedures, and therefore closed testing procedures, themselves. The difference between $\boldsymbol\phi$ and $\boldsymbol\psi$, if both are admissible, is that for every $S \subseteq I$, $\boldsymbol\phi_S$ uses only the local information in $\mathbf{d}^S_{\boldsymbol{\phi}}(S)$, but the same does not necessarily hold for $\boldsymbol\psi_S$.

In Section B in the Supplemental Information we give an example of a local improvement of an admissible monotone procedure. Local improvements are also possible in case null hypotheses are composite, using the Partitioning Principle, as shown in \cite{Finner2002}, examples 4.1--4.3, and \cite{Goeman2010}, section 4. For many well-known procedures, e.g.\ Holm's procedure under arbitrary dependence, we believe that local improvements do not exist. However, we have no general theory on the relationship between admissibility of a monotone procedure and admissibility of its local members. We leave this as an open problem.

\section{Consonance and familywise error} \label{sec consonance}

Theorem \ref{thm admissibility} establishes a necessary and sufficient condition for admissibility of monotone true discovery guarantee procedures, and therefore of FDP-controlling procedures. At first sight, our results may seem at odds with those of \cite{Romano2011}, who proved that for FWER control, which is a special case of the true discovery guarantee requirement, only consonant procedures are admissible. However, this seeming contradiction disappears when we realize that admissibility of a procedure as a true discovery guarantee procedure does not automatically imply admissibility as a FWER controlling procedure and vice versa. 

We call a procedure $\mathbf{d}^I$ \emph{consonant} if it has the property that for every $S \subseteq I$, $\mathbf{d}^I(S) > 0$ implies that  for at least one $i\in S$ we have $\mathbf{d}^I(\{i\}) = 1$, almost surely for all $\mathrm{P} \in \Omega$. Conceptually, consonant procedures allow pinpointing of effects. If $\mathbf{d}^I(S)>0$, signal has been detected somewhere in $S$. A consonant procedure in this case can always find at least one elementary hypothesis to pin the effect down on. This is a desirable property, as it can be unsatisfactory for a researcher to know that an effect exists but not where it can be found. However, \cite{Goeman2017} argued that for FDP control, non-consonant procedures can be far more powerful in large-scale multiple testing procedures than consonant ones.

In Section C of the Supplemental Information we go more deeply into the theory of consonant procedures in relation to admissibility of procedures as FWER controlling procedures. We extend the result of \cite{Romano2011}, showing that admissible FWER controlling procedures must be closed testing procedures with consonant local tests, but also closed testing procedures with admissible local tests. Conversely, if the local tests are both admissible and consonant, then the resulting closed testing procedure is admissible.

\section{Improving methods 1: \cite{Meinshausen2006a}} \label{sec HC}

Existing methods may be improved by embedding them in a closed testing procedure. We illustrate this with the method of Higher Criticism \citep{Donoho2004}, which defines a global test for the null hypothesis $H_I$, as follows. Let $I=\{1,\ldots,m\}$, and assume we have $p$-values $\mathbf{p}_1, \ldots, \mathbf{p}_m$, independent and stochastically larger than uniform under $H_I$. For this null hypothesis, Higher Criticism defines the test
\[ 
\boldsymbol\phi_I = \mathds{1}\Bigg\{ \max_{k_0 \leq j \leq k_1} \frac{\sqrt{m}(j/m - \mathbf{p}_{(j)})}{\sqrt{\mathbf{p}_{(j)}(1-\mathbf{p}_{(j)})}} \geq a_m \Bigg\},
\]
for suitably chosen $k_0$ and $k_1$, where $\mathbf{p}_{(1)} \leq \ldots \leq \mathbf{p}_{(m)}$ are the sorted $p$-values, and $a_m$ is a suitably chosen critical value. \cite{Donoho2004} proposed $a_m = (1+a)\sqrt{2\log\log(m)}$ for some $a>0$, assuming large $m$. Several finite-$m$ adjustments have been proposed \citep{Hall2010, Barnett2014}. We will use $k_0=1$ and $k_1=m$. \cite{Meinshausen2006a} improved upon Higher Criticism by showing that discoveries may also be counted, proving that
\[
\mathbf{f}_I = \Bigg\lceil \max_{t \in [0,1)} \frac{|\{i \in I\colon p_i \leq t\}| - mt - a_m \sqrt{mt(1-t) }}{1-t} \Bigg\rceil
\]
is a $(1-\alpha)$-lower confidence bound for the number of false hypotheses $|I_1|$. We have $\boldsymbol\phi_I = \mathds{1}\{\mathbf{f}_I > 0\}$, so $\mathbf{f}_I$ is consistent with the higher criticism test, and uniformly improves it as a true discovery guarantee procedure.

Can we improve $\mathbf{f}_I$ further? First, we can use (\ref{def interpolation}) to interpolate, getting
\begin{equation} \label{eq HC consonant}
\mathbf{d}^I(S) = \mathbf{f}_I - m + |S|.
\end{equation}
The resulting method is consonant, and by Corollary \ref{thm corol} it is equivalent to a closed testing procedure with local tests $\boldsymbol\psi_S = \mathds{1}\{|S| > m - \mathbf{f}_I \}$ for every $S \subseteq I$, where we note that $\boldsymbol\psi_I = \boldsymbol\phi_I$. The interpolated method improves upon $\mathbf{f}_I$ by giving non-trivial $\mathbf{d}^I(S)$ for $S \neq I$ with large $|S|$, but still has $\mathbf{d}^I(I) = \mathbf{f}_I$.

Further improvement is possible by noting that the suite $\boldsymbol\psi$ is not admissible. In fact, $\boldsymbol\psi$ is uniformly improved by $\boldsymbol\phi$, the suite of Higher Criticism local tests. This test is suggested by the recipe of Theorem \ref{thm sufficiency} for improving methods. In Section E of the Supplemental Information we show that
\begin{equation} \label{eq HC compare}
\boldsymbol\phi_S \geq \boldsymbol\psi_S \qquad \textrm{for all $S \subseteq I$},
\end{equation}
and that $\boldsymbol\phi_S$ uniformly improves $\boldsymbol\psi_S$ for $\emptyset \subset S \subset I$. It follows that $\mathbf{d}_{\boldsymbol\phi}^I$ uniformly improves  $\mathbf{d}^I$, and that even $\mathbf{d}_{\boldsymbol\phi}^I(I)$ uniformly improves $\mathbf{f}_I$ as a confidence bound for $|I_1|$, as we shall see.

To solve the issue of computing $\mathbf{d}_{\boldsymbol\phi}^I$, we write \citep{Gontscharuk2016}
\begin{equation} \label{eq simeslike}
\boldsymbol\phi_S = \mathds{1}\{\mathbf{p}_{(i\mathbin{:}S)} \leq l_{i\mathbin{:}|S|} \textrm{\ for at least one $i=1,\ldots, |S|$} \},
\end{equation}
where $\mathbf{p}_{(i\mathbin{:}S)}$, for $1 \leq i \leq |S|$, is the $i$th smallest $p$-value among the multiset $\{p_i\colon i\in S\}$, and
\begin{equation} \label{eq l HC}
l_{i\mathbin{:}s} = \frac{2i + a_s^2 - \sqrt{(2i + a_s^2)^2 - 4i^2 (s+a_s^2)/s}}{2(s+a_s^2)}.
\end{equation}
Written like this, we see that $\boldsymbol\phi$ is similar to the Simes tests investigated by \cite{Goeman2017}. For calculating $\mathbf{d}_{\boldsymbol\phi}^I(S)$ we can use a generalization of the algorithm presented in that paper, given as Lemma \ref{thm shortcut}.
\begin{lemma} \label{thm shortcut}
If $\boldsymbol\phi_S$, $\emptyset\neq S \subseteq I$, is of the form (\ref{eq simeslike}), with $l_{i\mathbin{:}m} \geq l_{i\mathbin{:}n}$ for all $i \geq 1$ and $0 \leq m \leq n$, then
\[
\boldsymbol\phi_S^I = \mathds{1}\{\mathbf{p}_{(i\mathbin{:}S)} \leq l_{i\mathbin{:}\mathbf{h}_I} \textrm{\ for at least one $i=1,\ldots, |S|$} \},
\]
and
\begin{equation} \label{eq CT KR}
\mathbf{d}_{\boldsymbol\phi}^I(S) = \max_{1\leq u \leq |S|} 1-u+|\{i\in S\colon \mathbf{p}_i \leq l_{u\mathbin{:}\mathbf{h}_I}\}|,
\end{equation}
where
\[
\mathbf{h}_I =\max \big\{s\in\{0,\ldots,|I|\}: \mathbf{p}_{(|I|-s+i\mathbin{:}I)}>l_{i\mathbin{:}s}, \textrm{ for } i=1,\ldots,s\big\}.
\]
\end{lemma}
The lemma offers calculation in quadratic time in the general case. For Higher Criticism, $\mathbf{h}_I$ can be calculated using bisection as in \cite{Goeman2017}, reducing computation time even to $O(m\log(m))$. We give a condition for the use of bisection in Section F of the Supplemental Information.

To illustrate the new method $\mathbf{d}_{\boldsymbol\phi}^I(S)$ we used a simple simulation using settings by \cite{Donoho2004}. We used $|I|= m=10^6$ independent one-sided $z$-tests. Of these, $10^3$ were under the alternative, with a mean shift of $\sqrt{0.30 \log(m)} \approx 2.04$. We used $a =1.08$ in the calculation of the critical value, which empirically gives good control of type I error for $m \approx 10^6$ and $\alpha=0.05$. We used $10^4$ replications. The power of Higher Criticsm in this setting is 98.0\%.

In this simulation, we found that $\mathbf{d}_{\boldsymbol\phi}^I(I)$ indeed improved Meinshausen and Rice's $\mathbf{f}_I$, although in this setting an improvement was found in only 2.2\% of the realizations. More importantly, however, the new $\mathbf{d}_{\boldsymbol\phi}^I(S)$ also makes meaningful statements for $S \neq I$.
\begin{figure}[!ht]
\includegraphics[width=\textwidth]{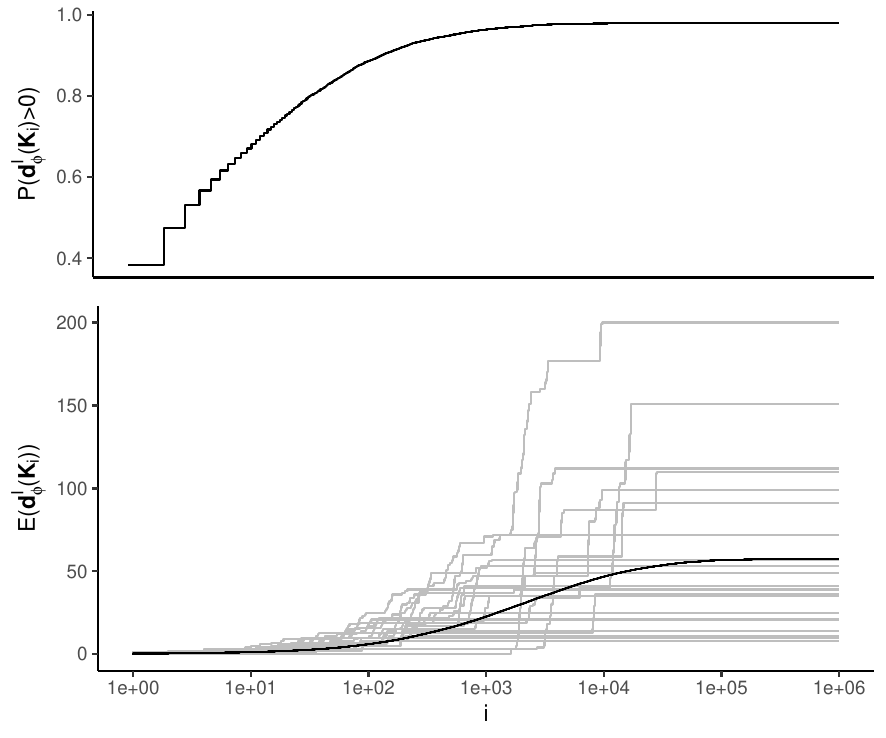}
\caption{Top: estimated probability of detection among the $i$ hypotheses with smallest $p$-values, using closed testing with the higher criticism local test for $m=10^6$ independent $p$-values. Bottom: Lower confidence bounds on the number of true discoveries among the same sets. Grey lines are 20 individual realizations; the black curve is the average over $10^4$ realizations.} \label{fig HC}
\end{figure}
Figure \ref{fig HC} (bottom) gives 20 realizations of $\mathbf{d}^I(\mathbf{K}_i)$ as a function of $i$, where $\mathbf{K}_i$ consists of the indices of the $i$ hypotheses with smallest $p$-values, with ties broken arbitrarily, as well as the expected curve. Figure \ref{fig HC} (top) gives the estimate of $\mathrm{P}(\mathbf{d}_{\boldsymbol\phi}^I(\mathbf{K}_i) > 0)$. We see that, by embedding Higher Criticism into a closed testing procedure, even with this weak signal we may make much stronger statements than only pure detection. In about 88.3\% of the realizations, we confidently detected signal within the set of 100 hypotheses with smallest $p$-values; in about 67.1\% this signal was in the top 10, and in 38.3\% of the realizations we even had a confident rejection of the single hypothesis with the smallest $p$-value. Substantial improvement of $\mathbf{d}_{\boldsymbol\phi}^I$ over $\mathbf{d}^I$, defined in (\ref{eq HC consonant}), is also clear since the latter gives $\mathbf{d}^I(S) = 0$ whenever $|S| \leq m-\mathbf{f}_I \approx 0.9999 \times 10^6$.

The strict line drawn by \cite{Donoho2004} between detectable and estimable effects is therefore, in our view, more like a gray zone. If we can detect effects, closed testing can count them. Though we may be unable to pinpoint effects, closed testing can close in on them.

\section{Improving methods 2: \cite{Katsevich2018}} \label{sec KR}

We chose as a second example a method recently proposed by \cite{Katsevich2018}. This elegant method (abbreviated K\&R) allows users to choose a $p$-value cutoff for significance post hoc, and uses stochastic process arguments to control both FDP and FDR. We focus on the FDP control property here. We use the same setting as in the previous section of $m$ $p$-values, independent and stochastically uniform under the null, and we use the same notation.

\citeauthor{Katsevich2018} showed that, if $\alpha \leq 0.31$,
\[
\mathrm{P}\big(|\mathbf{K}_i\cap I_0| \leq c(1+m\mathbf{p}_{(i)}) \textrm{\ for all $1\leq i \leq m$}\big) \geq 1-\alpha,
\]
where $c = -\log(\alpha)/\log(1-\log(\alpha))$, and $\mathbf{p}_{(i)}$ is the $i$th smallest $p$-value. For $\alpha = 0.05$ we have $c\approx 2.163$. As in Section \ref{sec procedures}, we can write this as a true discovery guarantee procedure on $I$ by writing
\begin{equation} \label{def KR}
\mathbf{d}^I(S) = \left\{ \begin{array}{ll} 0 \vee \lceil i - c(1+m\mathbf{p}_{(i)})\rceil & \textrm{if $S = \mathbf{K}_i$ for some $1\leq i \leq m$,} \\ 0 & \textrm{otherwise,} \end{array} \right.
\end{equation}
where we round up to ensure that $\mathbf{d}^I(S)$ is always an integer.

Is the procedure (\ref{def KR}) admissible, and if not, how can we improve it? We apply the results of this paper. First, we remark that the method as defined is not coherent. The interpolation of the procedure is given by
\begin{equation} \label{eq interpolation KR}
\mathbf{d}^I(S) = 0 \vee \max_{k=1,\ldots,|S|} \big\lceil k - c(1+m\mathbf{p}_{(k\mathbin{:}S)}) \big\rceil,
\end{equation}
taking $\mathbf{d}^I(\emptyset) = 0$ implicitly. The derivation of equation (\ref{eq interpolation KR}) is given in Section E of the Supplemental Information.

We note that interpolated method (\ref{eq interpolation KR}) makes non-trivial statements for sets $S$ not of the form $\mathbf{K}_i$, and may even improve $\mathbf{d}^I(\mathbf{K}_i)$ for some $i$. It may be checked using Lemma \ref{thm coherence} that the procedure (\ref{eq interpolation KR}) is coherent, so no further rounds of interpolation are needed. The K\&R procedure was not developed for a specific scale $m$. Writing $|I|$ for $m$ in (\ref{eq interpolation KR}) we have a procedure that is defined for general $I$, and it is easy to check that $(\mathbf{d}^I)_{I \subseteq \mathbf{\mathds{N}}, |I|<\infty}$ is monotone.

Next, we use Theorem \ref{thm sufficiency} to embed the method in a closed testing procedure, which results in further improvement of the procedure. By the theorem, the local test for finite $\emptyset \neq S \subseteq \mathbb{N}$ is given by
\[ 
\boldsymbol\phi_S = \mathds{1}\{\mathbf{p}_{(i\mathbin{:}S)} \leq (i-c)/c|S| \textrm{\ for at least one $i=1,\ldots, |S|$} \}.
\]
We will construct the closed testing procedure based on this local test. We note that the local test is of the form assumed in Lemma \ref{thm shortcut}, with
\begin{equation} \label{eq KR local}
l_{i\mathbin{:}s} = \frac{i-c}{cs}
\end{equation}
if $s \neq 0$, and $l_{i\mathbin{:}0} = 1$. By Theorem \ref{thm sufficiency}, the method (\ref{eq CT KR}) with (\ref{eq KR local}) is everywhere at least as powerful as the interpolated method (\ref{eq interpolation KR}). In fact, it is a uniform improvement of that method as we shall see in the simulation experiment below.

The next question is whether the method defined by (\ref{eq CT KR}) with (\ref{eq KR local}) is admissible, or whether it can be further improved. We can verify this using Theorem \ref{thm admissibility} by checking whether the local tests are admissible. It is immediately obvious that this is not the case. Taking e.g.\ $|S|=1$, we see that at $\alpha=0.05$ with $c\approx 2.163$ we have $\boldsymbol\phi_S = \mathds{1}\{\mathbf{p}_{(1\mathbin{:}S)} \leq (1-c)/c < 0\} = 0$, which is clearly not admissible. We may freely decrease $c$ to $c_1 = 1/(1+\alpha) \approx 0.952$ to obtain the uniformly more powerful local test $\boldsymbol\phi_S = \mathds{1}\{\mathbf{p}_{(1\mathbin{:}S)} \leq \alpha\}$. We can use the same reasoning for $|S|= 2,3,\ldots$, decreasing the value of $c$ to the minimal value that guarantees type I error control. This value may easily be calculated numerically since the worst case distribution of $(\mathbf{p}_i)_{i \in S}$ under $H_S$ is the independent uniform case. We obtain a new local test of the form (\ref{eq simeslike}) with
\begin{equation} \label{eq KR improved critical}
l_{i\mathbin{:}s} = \frac{i-c_s}{c_ss}.
\end{equation}
We tabulated the values of $c_s$ (taking $\alpha=0.05$) for some values of $s$ in Table \ref{tab cm}. Note that $c_s \leq c$ for all $s$. Since $l_{i\mathbin{:}s}$ is monotone in $c_s$ the new local test uniformly improves the old one. We note that with these choices of $c_s$ the critical values $l_{i\mathbin{:}s}$ cannot be further increased without destroying type I error control of the local tests, so we conclude that the resulting local tests are admissible, provided that the test $\mathds{1}\{\mathbf{p}_i \leq \alpha\}$ is admissible as an $\alpha$-level local test of $H_i$ for all $i$ and $\alpha$. Assuming this, by Theorem \ref{thm admissibility} the resulting true discovery guarantee procedure is admissible. We note that, since $c_s$ is increasing in $s$, (\ref{eq KR improved critical}) still fulfils the conditions of Lemma \ref{thm shortcut}, so that the admissible method is still computable using Lemma \ref{thm shortcut}.

\begin{table}[ht]
\caption{Values of $c_s$ calculated by Monte Carlo integration ($10^6$ samples)} \label{tab cm}
\centering
\begin{tabular}{rrrrrrrrrrrrrr}
\hline
$s$ & 1 & 2 & 3 & 4 & 5 & 7 & 10 & 15 & 20 & 50 & 100 & 500 & 1000\\
$c_s$   & 0.95 & 1.38 & 1.55 & 1.64 & 1.71 & 1.78 & 1.84 & 1.90 & 1.92 & 1.98 & 2.00 & 2.01 & 2.02 \\\hline
\end{tabular}
\end{table}

We have started with the procedure of \cite{Katsevich2018} and improved it uniformly in three steps: the method was first improved by interpolation. The resulting coherent method was further improved by embedding it in a closed testing procedure, and finally that closed testing procedure was improved to an admissible method by improving its local tests. This way we obtained a sequence of four methods, each uniformly improving the previous one. We will call them the \emph{original} (\ref{def KR}), \emph{coherent} (\ref{eq interpolation KR}), \emph{closed}, defined by (\ref{eq CT KR}) with (\ref{eq KR local}), and \emph{admissible} method, defined by (\ref{eq CT KR}) with (\ref{eq KR improved critical}). We performed a small simulation experiment to assess the relative improvement made with each of the three steps. We used $m=1000$ hypotheses, of which $m_0$ were true, and $m_1=m-m_0$ false. We sampled $p$-values independently. For true null hypotheses, we used $\mathbf{p}_i \sim \mathcal{U}(0,1)$. For false null hypotheses, we used $\mathbf{p}_i \sim \Phi^{-1}(-\gamma\mathbf{Z})$, where $\Phi$ is the standard normal distribution function, and $\mathbf{Z} \sim \mathcal{N}(0,1)$. We took values $m_1=8, 40, 200$ and $\gamma= 2, 2.5, 3$. A true discovery guarantee procedure gives exponentially many output values. We report only results for sets $\mathbf{K}_i$ of the $i$ smallest $p$-values, as the original method did. Calculation for the closed and admissible methods was in quadratic time based on Lemma \ref{thm shortcut}. We calculated $\mathbf{d}^I(\mathbf{K}_i)$, $i=1,\ldots,m$ for the closed and admissible methods in less than 0.1 seconds on a standard PC.

\begin{figure}[!htp]
\includegraphics[width=0.8\textwidth]{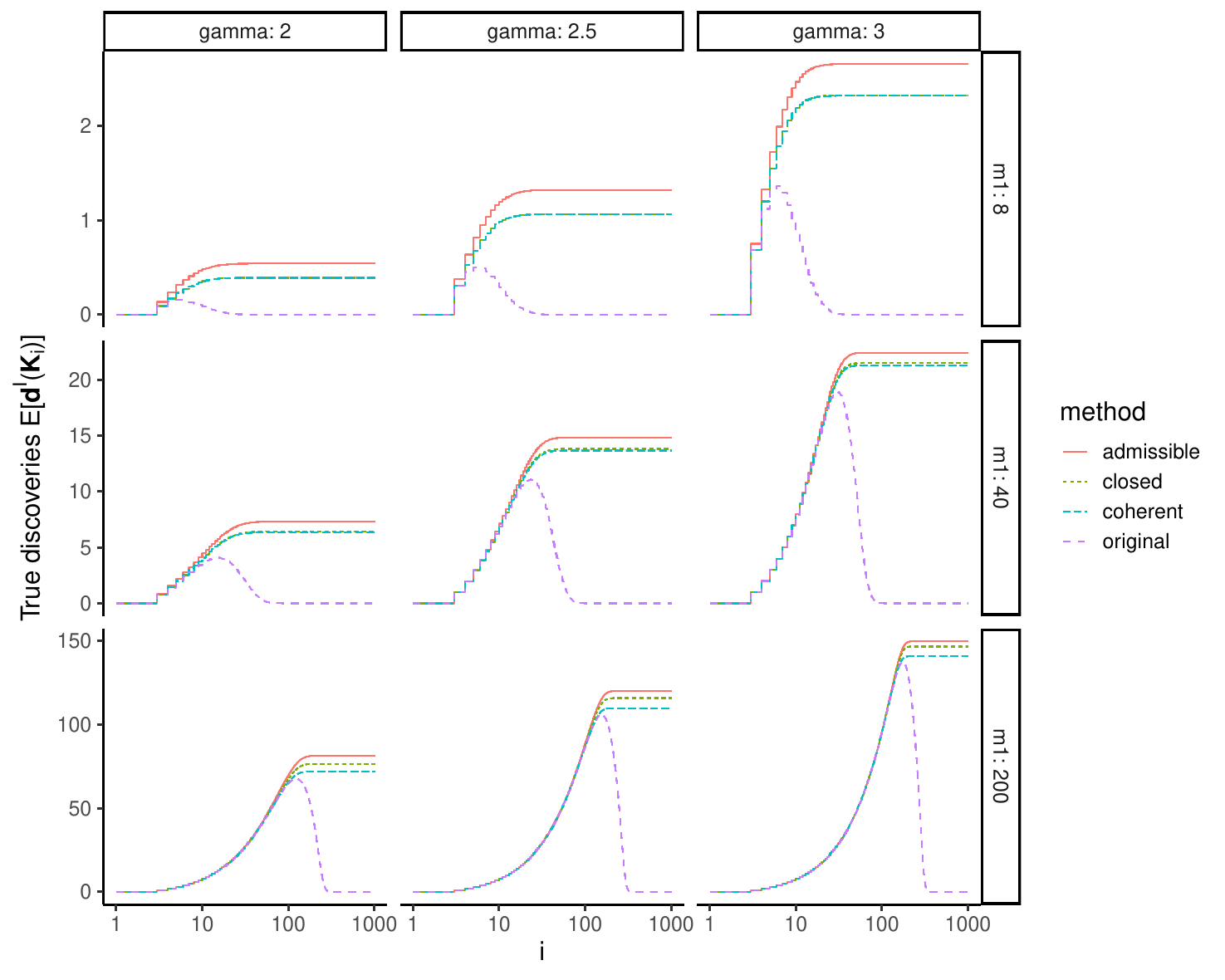}
\includegraphics[width=0.8\textwidth]{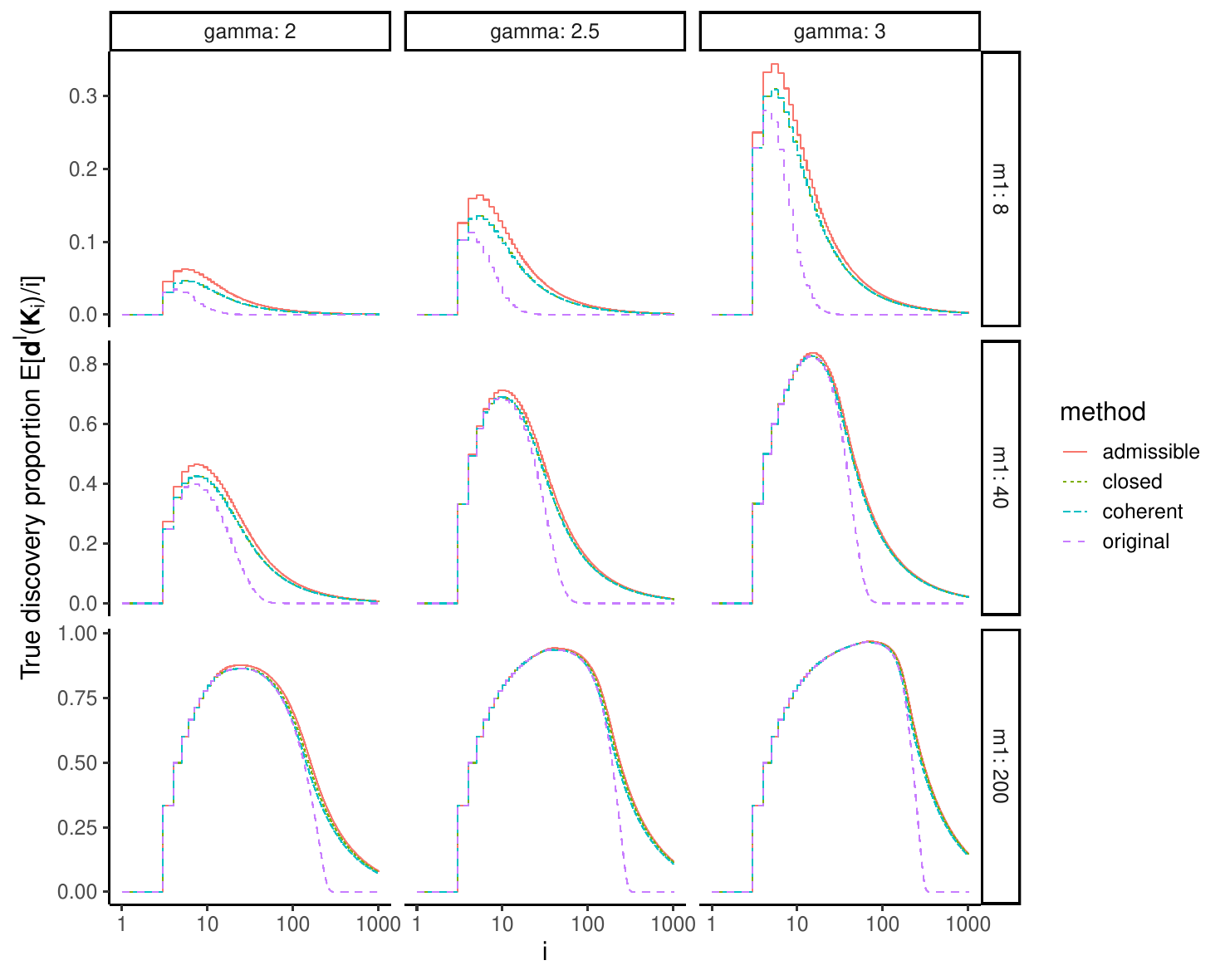}
\caption{Lower confidence bound on the number of true discoveries (top) and true discovery proportion (bottom) among the $i$ hypotheses with smallest $p$-values, relating to the method of \cite{Katsevich2018} and its successive uniform improvements. Curves are the averages over $10^4$ realizations. Since each method uniformly improves the previous one, any observed difference is automatically significant.} \label{fig KR}
\end{figure}

The results are given in Figure \ref{fig KR}, in terms of number of true discoveries $\mathbf{d}^I(\mathbf{K}_i)$ (top) and in terms of true discovery proportions $\mathbf{d}^I(\mathbf{K}_i)/i$ (TDP; bottom). For each setting and each method we report the average value of $\mathbf{d}^I(\mathbf{K}_i)/i$ over $10^4$ simulations. Several things can be noticed about these simulation results. 

The most important finding is that all three improvement steps can be substantial. The improvement from the original to the coherent procedure is perhaps largest. It is especially noticeable for large rejected sets, where the original method may all too often give $\mathbf{d}^I(S)=0$, especially if $|S|\gg m_1$. The peak of the TDP is improved if the TDP of the original method was low. The second improvement, from the coherent procedure to closed testing is substantial in terms of numbers of discoveries only if $m_1$ is large. This is natural because the improvement can be seen as a ``step-down'' argument, implicitly incorporating an estimate of $m_1$ into the procedure. Even with large $m_1$ it is the improvement is negligible on the TDP scale. 
The final improvement from the initial closed testing to the admissible procedure is clear throughout the figure. It is largest in terms of number of discoveries when $m_1$ is large, but largest in terms of TDP when $m_1$ is small.

Although the improvement from the coherent to the closed procedure seems the smallest one, we emphasize that closed testing was  crucial for the construction of the admissible procedure. We also note that the method of K\& R, as well as its improvements, make no useful FWER rejections at all: in Figure \ref{fig KR} we see that $\mathrm{E}[\mathbf{d}^I(\mathbf{K}_i)]\approx 0$ for $i \leq 2$. This phenomenon is analyzed in more detail in Section D of the Supplemental Information.

\section{Discussion}

We have studied the class of all methods controlling tail probabilities of false discovery proportions. This class encompasses very diverse methods, e.g.\ familywise error control procedures, false discovery exceedance procedures, joint error rate controlling methods, and cluster inference. We have shown that all such procedures can be written as methods simultaneously controlling false discovery proportions over all subsets of the family of hypotheses. This rewrite, trivial as it may be in some cases, is valuable in its own right, because it makes it possible to study methods jointly that seemed incomparable before, and takes a step in reducing the `plethora of error rates' lamented by \cite{Benjamini2010}. Moreover, methods that were constructed to give non-trivial error bounds for only a single random hypothesis set of interest, now give simultaneous error bounds for all such sets, allowing their use in flexible selective inference in the sense advocated by \cite{Goeman2011}.

We have formulated all such procedures in terms of a $(1-\alpha)$-true discovery guarantee, i.e.\ giving a $(1-\alpha)$-lower confidence bound to the number of true discoveries in each set, because this representation is mathematically easier to work with. Also, by emphasizing true rather than false discoveries, it gives a valuable positive frame to the multiple testing problem. Otherwise, this change in representation is purely cosmetic; we may continue to speak of FDP control procedures.

Admissibility is a very weak requirement for statistical tests, as under a weak assumption all tests that exhaust their $\alpha$-level are admissible. However, admissibility is not so easy to achieve for FDP control procedures. We have formulated a condition for admissibility of FDP control procedures that is both necessary and sufficient. All admissible FDP control procedures are closed testing procedures, and all closed testing procedures are admissible as FDP control procedures, provided they are well-designed in the sense that all their local tests are admissible. Apparently, control of false discovery proportions and closed testing procedures are so closely tied together that the relationship seems almost tautological. Admissibility is closely tied to optimality. Since optimal methods must be admissible, and admissible methods must be closed testing procedures, we have shown that only closed testing procedures can be optimal.

This theoretical insight has great practical value for methods designers. It can be used to uniformly improve existing methods, as we have demonstrated on the methods of \cite{Meinshausen2006a} and \cite{Katsevich2018}. Given a procedure that controls FDP, we first make sure it is coherent. Next, we can explicitly construct the local tests implied by the procedure, and turn it into a closed testing procedure. To check admissibility, we now only need to check admissibility of the local tests. Each step may result in substantive improvement, as we have shown in simulations. Alternatively, when designing a method we may start from a suite of local tests that has good power properties. The options are virtually unlimited here. The validity of the local test as an $\alpha$-level test guarantees control of FDP. Correlations between test statistics, that often complicate multiple testing procedures, should be properly taken into account by the local test. Admissibility of the local tests guarantees admissibility of the resulting procedure. In both cases the computational problem remains that closed testing may require exponentially many tests, but this is the only remaining problem. Polynomial time shortcuts are possible. Ideally these are exact, as for K\&R and higher criticism above, and admissibility is retained. If the full closed testing procedure is not computable for large testing problems, we may settle for an inadmissible but computable method, based on a conservative shortcut \citep[e.g.][]{Hemerik2018, Hemerik2018a}. It may still be worthwhile to compare such a method to full closed testing in small-scale problems to see how much power is lost.

Concretely, in Lemma \ref{thm shortcut} have given an exact computational shortcut that can be used for closed testing with a wide range of local tests, e.g.\ to the local tests implied by the False Discovery Rate controlling procedures of \cite{Blanchard2009}, to other local tests implied by the Dvoretzky-Kiefer-Wolfowitz inequality \citep{Genovese2004, Meinshausen2006a}, to the local tests implied by second and higher order generalized Simes constants \citep{Cai2008, Gou2014}, and to the local tests implied by the FDR controlling procedures of \cite{Benjamini1999}, and \citet[equation 4.1]{Romano2006}. Using the lemma, computation time of closed testing is quadratic, even reducing to linearithmic in some cases.

We have defined admissibility in terms of simultaneous FDP control for all possible subsets of the family of hypotheses. In some cases we may not be interested in all of these sets, as e.g.\ when targeting FWER control exclusively. Even when only interested in some of the subsets, we retain the result that admissible procedures must be closed testing procedures. We lose, however, the property that all such procedures are automatically admissible if they have admissible local tests. Additional criteria might come in, such as consonance in the case of familywise error control. Variants of consonance may be useful as well \citep{Brannath2010}.

Our focus has been mostly on monotone procedures. Such procedures are defined for multiple testing problems on different scales simultaneously. Connecting between different scales, they have the property that adding more hypotheses to the multiple testing problem will never result in stronger inference for the hypotheses that were already there. This is an intuitively desirable property by itself, which prevents some paradoxes \citep{Goeman2014}. Monotone procedures have additional valuable properties: viewed as closed testing procedures, they have local tests that are truly local: the local test on $S$ uses only the information that the corresponding local procedure $\mathbf{d}^S$ uses. Admissible monotone procedures, however, may sometimes be locally improved, and we have given an example of this. Such improvements, if admissible, must still be closed testing procedures with admissible local tests themselves.

We have restricted to finite testing problems. Extensions to countably infinite problems are of interest e.g.\ when considering online control \citep{Javanmard2018}. The results of this paper may trivially be extended to allow infinite $|I|$ if we are willing to assume that $|I_1|<\infty$, so that $\mathbf{d}^I<\infty$. If $|I_1|$ is unbounded, care must be taken to scale $\mathbf{d}$ properly to keep it in the non-trivial range. This scaling adds some technical complexity, and is not assumption-free because $\mathbf{d}^I(S)$ scales with the unknown $|S \cap I_1|$. However, since most of the results of this paper compare methods that obviously require the same scaling, we conjecture that the optimality of closed testing will translate to FDP control in countable and even uncountable multiple testing problems. We leave this to future research.

Finally, we remark that we have only considered procedures that control tail probabilities of the false discovery proportion. These methods can also be used for bounding the median FDP \citep{Goeman2011}. However, if there is interest in the central tendency of FDP it is more common to bound the mean FDP, better known as False Discovery Rate (FDR). Given the close connection we have established between closed testing and FDP tail probabilities, it is likely that there is  also a connection between closed testing and FDR control. Some connections have already been found between Simes-based closed testing and the procedure of \cite{Benjamini1995} by \cite{Goeman2017}. It is likely that there are more such connections. Any procedure that controls FDR, since FDR control implies weak FWER control, implies a local test and can therefore be used to construct a closed testing procedure. Conversely,  if FDP is controlled with $(1-\alpha)$-confidence at level $\gamma$, then FDR is controlled at $\alpha(1-\gamma)+\gamma$, as \cite{Lehmann2005} have shown. More profound relationships may be found in the future.

\appendix

\setcounter{theorem}{0}

\section{Existence of admissible procedures} \label{sec existence}

Admissibility as defined in Section \ref{sec monotone} is known in the literature as $\alpha$-admissibility on $\Omega$. Alternative definitions of admissibility exist \citep[Section 6.7]{Lehmann2006}. With $\alpha$-admissibility on $\Omega_1 = \Omega \setminus H$, the law $\mathrm{P}$ under which $\tilde{\boldsymbol{\delta}}$ improves $\boldsymbol{\delta}$ with positive probability must be in $\Omega_1$. With $d$-admissibility, there is an additional requirement that $\tilde{\boldsymbol{\delta}}$ may not improve $\boldsymbol{\delta}$ with positive probability for $\mathrm{P} \in H$. However, $\alpha$-admissibility on $\Omega$ is most commonly considered in the multiple testing context \citep[e.g.][Section 9.3]{Lehmann2006}, because with multiple hypotheses there is no unique $\Omega_1$.

Admissible tests may not always exist. Consider for example the model where $\mathbf{X} \sim \mathrm{P}_\theta = \mathcal{N}(\mu, \sigma^2)$, with parameter space $\Theta = (-\infty, \infty) \times [0, 1]$ for $\theta = (\mu, \sigma^2)$. Let $\boldsymbol{\delta}$ be any valid test for $H: \mu=0$, for example $\boldsymbol\delta = \mathds{1}\{\mathbf{X} > z_{1-\alpha}\}$, where $z_{1-\alpha}$ is the $1-\alpha$-quantile of the standard normal distribution. Then for any $c$ outside the rejection region of $\boldsymbol\delta$,
\[
\tilde{\boldsymbol{\delta}} = \boldsymbol{\delta} + \mathds{1}\{\mathbf{X} = c\}
\]
improves $\boldsymbol{\delta}$ with positive probability for $\mu=c$ and $\sigma^2=0$. Since any valid test may be improved in a similar way unless $\alpha=1$, an admissible test does not exist.

However, we can easily guarantee existence of admissible tests if we rule out degenerate models. A weak assumption for this is the following.

\begin{assumption}[common null events] \label{ass shared null}
For every measurable event $E$ and for every $\mathrm{P}_1, \mathrm{P}_2 \in \Omega$ we have that $\mathrm{P}_1(E) = 0$ implies that $\mathrm{P}_2(E) = 0$.
\end{assumption}

Under Assumption \ref{ass shared null} the collection of null events is common to all parameter values; there are no events that happen with positive probability for some parameter values, but with probability zero for others. This is a weak assumption that holds for most models regularly used in applied statistics, both continuous and discrete, when we are willing to exclude deterministic corner cases from the parameter space. In the example above, Assumption \ref{ass shared null} holds if we simply restrict the parameter space to $\Theta = (-\infty, \infty) \times (0, 1]$, excluding $\sigma^2=0$.

If we accept Assumption \ref{ass shared null}, Lemma \ref{thm adm test} presents a very simple sufficient condition for admissibility: every statistical test of $H$ that fully exhausts its $\alpha$-level for some $\mathrm{P}\in H$ is admissible.

\begin{lemma} 
If Assumption \ref{ass shared null} holds, then a statistical test $\boldsymbol\delta$ of hypothesis $H \subseteq \Omega$ is admissible if $\mathrm{P} \in H$ exists such that $\mathrm{P}(\boldsymbol\delta=1) = \alpha$.
\end{lemma}

\section{A local improvement} \label{sec local improvement}

In this section we construct a local improvement of an admissible monotone procedure to illustrate Proposition \ref{thm locally admissible}. Assume that for each $H_i$, $i \in C$ we have a $p$-value $\mathbf{p}_i$. Assume that each $\mathbf{p}_i$ is standard uniform if $H_i$ is true. Under these assumptions we can define the standard fixed sequence testing procedure, which starts testing $H_1$ using $\mathbf{p}_1$ at level $\alpha$, continues one by one with $H_2$, $H_3$ in order, and stops when it fails to reject some hypothesis. It is well known that this procedure is a closed testing procedure. The local test is defined for all $S \subseteq C$ by
\[
\boldsymbol\phi_{S} = \mathds{1}\{\mathbf{p}_{\min(S)} \leq \alpha\}.
\]
If we assume that the test $\boldsymbol\phi_{\{i\}} = \mathds{1}\{\mathbf{p}_i \leq \alpha\}$ is admissible for all $i$, then, by Theorem \ref{thm admissibility}, the fixed sequence procedure is admissible.

We will now make some additional assumptions that will allow a uniform improvement of the procedure at the fixed scale $I\subseteq C$. Assume for convenience that all $p$-values are independent. Next, assume that the distribution of every $\mathbf{p}_i$ is constrained even under the alternative. Assume that some $h(x)$ exists such that, for every $\mathrm{P} \in \Omega$,
\[
\mathrm{P}(\mathbf{p}_i \leq x) \leq h(x).
\]
This means that the power for each test is inherently limited. Even under the alternative, we reject e.g.\ $H_1$ with probability at most $h(\alpha)$. Clearly $h(x) \geq x$ or we would not have uniformity under the null. We will now demonstrate that the fixed sequence procedure can be uniformly improved, locally at any $I$ with $|I|>1$, if $h(\alpha)<1$.

We start with the simple case $I = \{1,2\}$. By Proposition \ref{thm locally admissible}, the improvement must be a closed testing procedure that involves local tests $\boldsymbol\psi_S \geq \boldsymbol\phi_{S}^I$. Consider $S = \{2\}$. Then $\boldsymbol\phi_{S}^I = \mathds{1}\{\mathbf{p}_1\leq\alpha, \mathbf{p}_2\leq\alpha\}$. Under $H_2$ this has $\mathrm{P}(\boldsymbol\phi_{S}^I =1) \leq \alpha h(\alpha)< \alpha$. Clearly, there is room for improvement. Let us consider the procedure with $\boldsymbol\psi_S =  \boldsymbol\phi_S$, for all $S \subseteq C$ except $S = \{2\}$, when
\[
\boldsymbol\psi_{S} = \mathds{1}\{\mathbf{p}_1\leq\alpha, \mathbf{p}_2\leq\alpha/h(\alpha)\}.
\]
The latter is a valid local test of $H_2$. The resulting procedure at $I=\{1,2\}$ starts testing $H_1$ at level $\alpha$, and continues, if $H_1$ is rejected, to test $H_2$ at level $\alpha/h(\alpha) > \alpha$. This is clearly a uniform improvement of the original procedure at $I=\{1,2\}$. To see that this is not a counterexample to Theorem \ref{thm admissibility}, consider $I=\{2\}$ instead. Clearly, we do not have
\[
\mathds{1}\{\mathbf{p}_1\leq\alpha, \mathbf{p}_2\leq\alpha/h(\alpha)\}\ =\ \mathbf{d}^{\{2\}}_{\boldsymbol\psi}(\{2\})\ \geq\ \mathbf{d}^{\{2\}}_{\boldsymbol\phi}(\{2\})\ =\
\mathds{1}\{\mathbf{p}_2\leq\alpha\}.
\]
The local uniform improvement at $I=\{1,2\}$ comes at the cost of a potential deterioration at $I=\{2\}$.

Similar local improvements actually exist for every finite $I \subseteq C$ with $|I| > 1$. Define recursively
\[
\alpha_1 = \alpha;\quad  \alpha_{i+1} = \alpha_i/h(\alpha_i).
\]
From this, fix some $I \subseteq C$, and define a local test as
\[
\boldsymbol\psi_S = \mathds{1}\{\mathbf{p}_i \leq \alpha_i \textrm{\ for all $i \in L_S$}\},
\]
where $L_S = \{i \in I\colon i \leq \min(S)\}$. To check that this is a valid local test, we verify that for all $\mathrm{P} \in H_S$
\[
\mathrm{P}(\boldsymbol\psi_S = 1)\ \leq\ \alpha_l\prod_{i=1}^{l-1} h(\alpha_i)\ =\ \alpha,
\]
where $l = |L_S|$. The resulting procedure is still a fixed sequence procedure that tests all $H_i$, $i\in I$ in order, stopping the first time it fails to reject. Only, rather than testing at level $\alpha$ every time, it tests at level $\alpha_i$ in step $i$. If $x < h(x) < 1$ for all $0 < x < 1$ the sequence $\alpha_1, \alpha_2, \ldots$ is strictly increasing and approaches 1.

Crucial for this example is the assumption that we have limited power and, more importantly, that we know the limit to the power. If we are not willing to assume that $h(\alpha)<1$, or if we do not know $h$, then the above local improvements are not possible. It is difficult to think of uniform local improvements in the case $h(x)=1$, and we believe they do not exist. It may be worthwhile to think of adaptive procedures that learn $h$ as the procedure moves along, but we will not pursue this direction here. In any case, due to the cost inherent to learning $h$, such a procedure would not uniformly improve $\mathbf{d}_{\boldsymbol\phi}$.

\section{Admissibility of FWER controlling procedures} \label{sec FWER}

In this section we take a sidestep to FWER control, investigating the concept of consonance, and extending some of the results of \cite{Romano2011} on admissibility of FWER controlling procedures. Consonance was defined in Section \ref{sec consonance}: we call a procedure $\mathbf{d}^I$ \emph{consonant} if it has the property that for every $S \subseteq I$, $\mathbf{d}^I(S) > 0$ implies that  for at least one $i\in S$ we have $\mathbf{d}^I(\{i\}) = 1$, almost surely for all $\mathrm{P} \in \Omega$. If $\mathbf{d}^I = \mathbf{d}_{\boldsymbol\phi}^I$, this definition is equivalent to the more usual formulation in terms of the suite $\boldsymbol\phi$, that $\boldsymbol\phi_{S}^I = 1$ implies that for at least one $i \in S$ we have $\boldsymbol\phi_{\{i\}}^I = 1$, almost surely for all $\mathrm{P} \in \Omega$. We call a monotone procedure $\mathbf{d}$ consonant if all local members $\mathbf{d}^I$, $I \subseteq C$ finite, are consonant. We call a suite $\boldsymbol\phi$ consonant if all $\boldsymbol\phi_S$, $S \subseteq C$ finite, are consonant.

For consonant procedures a stronger version of Lemma \ref{thm coherence} holds.

\begin{lemma} \label{thm VW consonant}
$\mathbf{d}^I$ is consonant and coherent if and only if, for every disjoint $V,W \subseteq I$,
\[ 
\mathbf{d}^I(V \cup W) = \mathbf{d}^I(V) + \mathbf{d}^I(W).
\]
\end{lemma}

Classically, focus in the literature on closed testing has been on FWER controlling procedures \citep{Henning2015}. A FWER-controlling procedure on a finite $I \subseteq C$ returns a set $\mathbf{R}^I \subseteq I$ such that, for all $\mathrm{P} \in \Omega$,
\[
\mathrm{P}(|\mathbf{R}^I \cap I_0| = 0) \geq 1-\alpha.
\]
As argued in Section \ref{sec procedures}, we can relate FWER controlling procedures to true discovery guarantee procedures and vice versa. If $\mathbf{R}^I$ is a FWER controlling procedure, then $\mathbf{r}^I$ with
\[
\mathbf{r}^I(S) = |S \cap \mathbf{R}^I|,
\]
for all $S\subseteq I$, is a coherent true discovery guarantee procedure, as we know from (\ref{eq interpolate FWER}). Conversely, if $\mathbf{d}^I$ is a coherent true discovery guarantee procedure, then
\[
\mathbf{R}^I =
\{i \in I \colon \mathbf{d}^I(\{i\}) = 1\}
\]
is a FWER controlling procedure. Both types of procedures may be created from local tests. The FWER controlling procedure from the suite $\boldsymbol\phi$ is given by
\begin{equation} \label{eq closed FWER}
\mathbf{R}^I_{\boldsymbol\phi} =
\{i \in I \colon \boldsymbol\phi^I_{\{i\}} = 1\} =
\{i \in I \colon \mathbf{d}_{\boldsymbol\phi}^I(\{i\}) = 1\}.
\end{equation}

We can compare the procedure $\mathbf{d}_{\boldsymbol\phi}$ defined from $\boldsymbol\phi$ through (\ref{def d}) with the procedure
\[
\mathbf{r}^I_{\boldsymbol\phi}(S) = |S \cap \mathbf{R}_{\boldsymbol\phi}^I|,
\]
indirectly defined through $\mathbf{R}_{\boldsymbol\phi}^I$. This is the procedure that discards all information in $\mathbf{d}^I_{\boldsymbol\phi}$ that is not contained in $\mathbf{R}^I_{\boldsymbol\phi}$. Lemma \ref{thm r inadmissible} describes consonance as the property that no information is lost in the process.

\begin{lemma} \label{thm r inadmissible}
If $\mathbf{d}^I_{\boldsymbol\phi}$ is consonant, $\mathbf{d}^I_{\boldsymbol\phi} = \mathbf{r}^I_{\boldsymbol\phi}$; otherwise $\mathbf{d}^I_{\boldsymbol\phi}$ uniformly improves $\mathbf{r}^I_{\boldsymbol\phi}$.
\end{lemma}

If FWER control is what we are after, however, we must look at admissibility of $\mathbf{R}^I$ directly. As with true discovery guarantee procedures, we will focus on monotone (stacks of) procedures defined for all finite $I \subseteq C$. We will call a procedure $\mathbf{R} = (\mathbf{R}^I)_{I \subseteq C, |I|<\infty}$ monotone if for all finite $J \subseteq I \subseteq C$ we have
\[
\mathbf{R}^J \supseteq \mathbf{R}^I \cap J.
\]
As above for true discovery guarantee procedures, it asserts that enlarging the multiple testing problem from $J$ to $I$ will never increase the number of rejections in $J$ \citep{Bretz2009, Goeman2010}. Analogous to the definition in Section \ref{sec monotone}, we define a uniform improvement of a monotone FWER control procedure $\mathbf{R}$ as a monotone FWER control procedure $\mathbf{\tilde R}$ such that (1.) $\mathbf{\tilde R}^I \supseteq \mathbf{R}^I$ for all finite $I \subseteq C$; and (2.) $\mathrm{P}(\mathbf{\tilde R}^I \supset \mathbf{R}^I) > 0$ for some $\mathrm{P} \in \Omega$ and some finite $I \subseteq C$. A procedure is admissible if no uniform improvement exists. What can we say about admissibility of FWER control procedures?

\cite{Romano2011} showed that consonance is necessary for admissibility of FWER controlling procedures. Proposition \ref{thm Romano} is a variant of their result for monotone procedures.

\begin{prop} \label{thm Romano}
If $\mathbf{R}$ is monotone and admissible, then a consonant suite $\boldsymbol\phi$ exists such that $\mathbf{R} = \mathbf{R}_{\boldsymbol\phi}$.
\end{prop}

We also have a second necessary condition for admissibility.

\begin{prop} \label{thm fwer admissible}
If $\mathbf{R}$ is monotone and admissible, then an admissible suite $\boldsymbol\psi$ exists such that $\mathbf{R} = \mathbf{R}_{\boldsymbol\psi}$.
\end{prop}

It would be tempting to conclude from Propositions \ref{thm Romano} and \ref{thm fwer admissible} that if $\mathbf{R}$ is admissible, then $\mathbf{R} = \mathbf{R}_{\boldsymbol\phi}$ with $\boldsymbol\phi$ consonant and admissible. Certainly under weak assumptions we may choose $\boldsymbol\phi = \boldsymbol\psi$. For example, if some $\boldsymbol\phi_S$ is inadmissible because it fails to exhaust the $\alpha$-level we may choose $\boldsymbol\psi_S$ as $\boldsymbol\phi_S$ plus a randomized multiple of $\boldsymbol\phi_{\{i\}}$ for some $i \in S$, if we allow randomized tests, and use the lemma from Section A in the Supplemental Information. However, we were unable to prove in full generality that $\boldsymbol\phi = \boldsymbol\psi$ is always possible. Perhaps in some awkward models admissible local tests cannot be consonant, and consonant tests cannot be admissible. We leave the question open when $\boldsymbol\phi = \boldsymbol\psi$ is possible. In converse however, if we can find a $\boldsymbol\phi$ that is both admissible and consonant, we have an admissible procedure:

\begin{prop} \label{thm cons and adm}
If $\boldsymbol\phi$ is consonant and admissible, then $\mathbf{R}_{\boldsymbol\phi}$ is admissible as a monotone procedure.
\end{prop}

\section{Some properties of the method of \cite{Katsevich2018}}

We may derive some the properties of the K\&R method and its improvements. We have that $\mathbf{d}^I(S) \leq |S|-2$ with probability 1 for the original method (\ref{def KR}), since $c >2$ if $\alpha=0.05$. This also holds for the coherent method (\ref{eq interpolation KR}). For the closed method the same is not true, but we have that if $\max_{1\leq i \leq m} \mathbf{p}_i > l_{3\mathbin{:}3} \approx 0.004$ we have $\boldsymbol\phi_I(S)=0$ for every $S$ with $|S|=2$, since $l_{1\mathbin{:}3}<l_{2\mathbin{:}3}<0$. Therefore
\[
\mathrm{P}(\mathbf{d}^I(S) \leq |S|-2 \textrm{\ for all $S\subseteq I$})\ \geq\ \mathrm{P}(\max_{1\leq i \leq m} \mathbf{p}_i > l_{3\mathbin{:}3})\ \approx\ 1,
\]
unless all hypotheses are false. For the admissible method, by an analogous reasoning using Lemma \ref{thm shortcut}, the same holds if $c_{|I_0|} > 2$, since then $c_{\mathbf{h}_I} > 2$ with large probability. This happens from $|I_0| \approx 100$. It follows that none of the methods in this section, not even the admissible method, should be expected make any FWER-rejections in practical applications. The admissible method is (almost) fully non-consonant in the sense that for all $i\in I$, $\boldsymbol\phi_{\{i\}}^I\approx 0$, and we have $\mathbf{R}_{\boldsymbol\phi}^I = \emptyset$ with probability almost 1 unless $|I_1|\approx|I|$. By \cite{Romano2011} the method is clearly inadmissible as a FWER-controlling method. By Theorem \ref{thm admissibility} it is admissible, however, as a true discovery guarantee procedure: its lack of power for FWER-type statements is compensated by larger power for non-FWER-type statements. Indeed, \citeauthor{Katsevich2018} have shown that their method may significantly outperform Simes-based closed testing \citep{Goeman2017} in some scenarios, which in turn outperforms consonant FWER-based testing in terms of FDP in large-scale testing problems.

\setcounter{lemma}{0}
\setcounter{prop}{0}

\section{Proofs of the theorems, propositions and lemmas} \label{sec proofs}

\begin{lemma} 
$\mathbf{g}^I_{\phi} = \mathbf{d}^I_{\phi}$.
\end{lemma}

\begin{proof}
Take any $S \in 2^I$. For any $V$ with $\boldsymbol{\phi}_V = 0$ there exists a $U = S \cap V \subseteq S$ which has $\boldsymbol{\phi}^I_U \leq \boldsymbol\phi_V = 0$ and $|S \setminus U| = |S \setminus V|$. Consequently, $\mathbf{d}^I_{\phi}(S) \leq \mathbf{g}^I_{\phi}(S)$.

For any $U \subseteq S$ with $\boldsymbol{\phi}^I_U = 0$ there is a $V \supseteq U$ with $\boldsymbol{\phi}_V = 0$ that has $|S \setminus V| \leq |S \setminus U|$. Consequently, $\mathbf{g}^I_{\phi}(S) \leq \mathbf{d}^I_{\phi}(S)$.
\end{proof}

\begin{lemma} 
If $\mathbf{d}^I$ is a true discovery guarantee procedure then so is $\mathbf{\bar d}^I$.
\end{lemma}

\begin{proof}
Let $E$ be the event that $\mathbf{d}^I(S) \leq |S \cap I_1|$ for all $S \subseteq I$. Suppose that $E$ happened and choose any $S \subseteq I$ and $U \subseteq I$. Then
\[
 \mathbf{d}^I(U) - |U \setminus S| + \mathbf{d}^I(S \setminus U) \leq |U \cap I_1| - |U \setminus S| + |(S \setminus U)\cap I_1| \leq |S \cap I_1|.
\]
Consequently, if $E$ happened, $\mathbf{\bar d}^I(S) \leq |S \cap I_1|$ for all $S \subseteq I$. Since $\mathrm{P}(E) \geq 1-\alpha$, we have coverage for the true discovery guarantee procedure $\mathbf{\bar d}^I$.
\end{proof}

\begin{lemma} 
$\mathbf{d}^I$ is coherent if and only if for every disjoint $V,W \subseteq I$ we have
\begin{equation} \label{eq coh disjoint}
\mathbf{d}^I(V) + \mathbf{d}^I(W) \leq \mathbf{d}^I(V \cup W) \leq \mathbf{d}^I(V) + |W|.
\end{equation}
\end{lemma}

\begin{proof} Suppose $\mathbf{d}^I$ is coherent and let $V$, $W$ be disjoint. Then, taking $S=V \cup W$ and $U=V$ in (\ref{def_coh}), we have
 \[
\mathbf{d}^I(V \cup W) \geq \mathbf{d}^I(V) - |V\setminus S| + \mathbf{d}^I(S \setminus V) = \mathbf{d}^I(V) + \mathbf{d}^I(W).
\]
Also, taking $S = V$ and $U=V\cup W$ in (\ref{def_coh}), we obtain $\mathbf{d}^I(V) \geq \mathbf{d}^I(V \cup W) - |W|$.

Next, suppose (\ref{eq coh disjoint}) holds for all disjoint $V, W$. For every $S \subseteq I$, we have, by (\ref{def interpolation}) and the left-hand inequality of (\ref{eq coh disjoint}),
\[
\mathbf{\bar d}^I(S) \leq \max_{U \in 2^I} \Big\{ \mathbf{d}^I(U \cup S) - |U \setminus S| \Big\}.
\]
Since $\mathbf{d}^I(U \cup S) \leq \mathbf{d}^I(S) + |U \setminus S|$ by the right-hand inequality of (\ref{eq coh disjoint}), we get $\mathbf{\bar d}^I(S) \leq \mathbf{d}^I(S)$. Since also $\mathbf{\bar d}^I(S) \geq \mathbf{d}^I(S)$, by taking $U=S$ in (\ref{def interpolation}), we have $\mathbf{\bar d}^I(S) = \mathbf{d}^I(S)$ for all $S \subseteq I$, so $\mathbf{d}^I$ is coherent.
\end{proof}

\begin{lemma} 
The procedure $\mathbf{d}_\phi^I$ is coherent.
\end{lemma}

\begin{proof} We use Lemma \ref{thm coherence}. Let $V,W \subseteq I$ be disjoint. Then some $U \subseteq V \cup W$ exists such that $\phi_U^I = 0$ and
\[
|V\setminus U|+ |W\setminus U| = |(V \cup W)\setminus U| = \mathbf{d}_\phi^I(V \cup W).
\]
Since $\phi_{U\cap V}^I \leq \phi_{U}^I = 0$, we have $\mathbf{d}_\phi^I(V) \leq |V\setminus U|$. Similarly, $\mathbf{d}_\phi^I(W) \leq |W\setminus U|$, so we have $\mathbf{d}_\phi^I(V)+\mathbf{d}_\phi^I(W)\leq \mathbf{d}_\phi^I(V\cup W)$.

Also, there exists $T \subseteq V \subseteq (V \cup W)$ such that $\boldsymbol\phi^I_T = 0$ and $|V\setminus T| = \mathbf{d}_\phi^I(V)$. Now
\[
|(V \cup W) \setminus T| = |V\setminus T| + |W| = \mathbf{d}_\phi^I(V) + |W|,
\]
so we have $\mathbf{d}_\phi^I(V\cup W) \leq \mathbf{d}_\phi^I(V) + |W|$.
\end{proof}

\begin{lemma} 
The procedure $\mathbf{d}_{\boldsymbol\phi} = (\mathbf{d}_{\boldsymbol\phi}^I)_{I \subseteq C, |I|<\infty}$ is a monotone procedure.
\end{lemma}

\begin{proof}
\cite{Genovese2006} already proved that $\mathbf{d}_\phi^I$ is a true discovery guarantee procedure for all finite $I \subseteq C$, and we have coherence by Lemma \ref{thm CT coherent}, so we only need to prove monotonicity. This is trivial from (\ref{def g}) and Lemma \ref{thm GW}. Take any finite $S \subseteq I \subseteq J \subseteq C$. Then
\[
\mathbf{d}^I_{\boldsymbol\phi}(S) = \min_{V \in 2^I} \{|S \setminus V|\colon \boldsymbol\phi_V =0\} \geq \min_{V \in 2^J} \{|S \setminus V|\colon \boldsymbol\phi_V =0\} = \mathbf{d}^J_{\boldsymbol\phi}(S).
\]
\end{proof}

\begin{theorem} 
Let $\mathbf{d}$ be a monotone procedure. Then, for every finite $S\subseteq C$,
\[
\boldsymbol\phi_S = \mathds{1}\{\mathbf{d}^{S}(S) > 0\}
\]
is a valid local test of $H_S$. For the suite $\boldsymbol\phi = (\boldsymbol\phi_S)_{S\subseteq C, |S|<\infty}$ we have, for all $S \subseteq I\subseteq C$ with $|I|<\infty$,
\[
\mathbf{d}^I_{\boldsymbol\phi}(S) \geq \mathbf{d}^I(S).
\]
\end{theorem}

\begin{proof} Take any finite $S \subseteq C$. Since $\mathbf{d}^{S}$ is a true discovery guarantee procedure on $S$, we have
\[
\max_{\mathrm{P} \in H_S} \mathrm{P}(\boldsymbol\phi_S = 1) = \max_{\mathrm{P} \in H_S} \mathrm{P}(\mathbf{d}^{S}(S) >  0) = \max_{\mathrm{P} \in H_S} \mathrm{P}(\mathbf{d}^{S}(S) >  |S \cap I_1|) \leq \alpha,
\]
so $\boldsymbol\phi_S$ is a valid test of $H_S$. This proves the first statement.

Take any finite $S \subseteq C$ again. We have $\boldsymbol{\phi}_S^I = 1$ if and only if $\mathbf{d}^W(W) > 0$ for all $S \subseteq W \subseteq I$. For all such $W $ we have, by coherence of $\mathbf{d}^W$ and monotonicity of $\mathbf{d}$,
\[
\mathbf{d}^W(W) \geq \mathbf{d}^W(S) \geq \mathbf{d}^I(S).
\]
Consequently, $\boldsymbol{\phi}^I_S \geq \mathds1\{\mathbf{d}^I(S) > 0\}$.
We obtain
\begin{equation} \label{eq d proof suf}
\mathbf{d}_{\boldsymbol \phi}^I(S) \geq \min_{U \in 2^S} \{|S \setminus U|\colon \mathbf{d}^I(U) = 0\}.
\end{equation}
By coherence of $\mathbf{d}^I$ and Lemma \ref{thm coherence}, we have, for all $U \subseteq S\subseteq I$, that
\[
|S \setminus U| \geq \mathbf{d}^I(S) - \mathbf{d}^I(U).
\]
Combining this with (\ref{eq d proof suf}) the second statement of the theorem follows.
\end{proof}

\begin{theorem} 
For every closed testing procedure there exists a partitioning procedure that rejects exactly the same hypotheses. For every partitioning procedure there exists a closed testing procedure that rejects exactly the same hypotheses.
\end{theorem}

Before we prove this theorem, we must first define the general partition procedure, following \citet{Finner2002}. Letting the hypotheses of interest be $(H_i)_{i \in I}$ as usual, define for every $S \subseteq I$ the partitioning hypothesis 
\[
\tilde H_S = H_S \setminus \bigcup_{i \in I \setminus S} H_i.
\]
This hypothesis is true if and only if all $H_i$, $i \in S$ are true and all $H_j$, $j \in I \setminus S$ are false. If, for every $S \subseteq I$, we have a valid statistical test $\boldsymbol{\psi}_S$ for $\tilde H_S$, then the partitioning procedure rejects $H_S$ if and only if 
\[
\tilde{\boldsymbol{\psi}}_S = \min\{\boldsymbol{\psi}_U\colon S \subseteq U \subseteq I\} = 1.
\]
The proof of validity of partitioning is similar to that of closed testing: let $T = \{i \in I\colon \mathrm{P} \in H_i\}$, then $\tilde H_T$ is true, and if $\boldsymbol{\psi}_T = 0$, which happens with probability at least $1-\alpha$, then no true hypothesis is rejected.

Now we can prove Theorem \ref{thm equiv_prin}

\begin{proof}
Consider the closed testing procedure defined by the suite $(\boldsymbol{\phi}_S)_{S \subseteq I}$. Define $\boldsymbol{\psi}_U = \boldsymbol{\phi}_U$ for all $U \subseteq I$. Since for all $U\subseteq I$, $\tilde{H}_U \subseteq H_U$, $\boldsymbol{\psi}_U$ is a valid test for $\tilde{H}_U$, we may define a partitioning procedure from $(\boldsymbol{\psi}_S)_{S \subseteq I}$. Since have
\[
\tilde{\boldsymbol{\psi}}_S = \min\{\boldsymbol{\phi}_U\colon S \subseteq U \subseteq I\} = \boldsymbol{\phi}_S^I,
\]
this partitioning procedure rejects exactly the same hypotheses as the closed testing procedure. This proves the first statement.

Consider a partitioning procedure defined by tests $(\boldsymbol{\psi}_S)_{S \subseteq I}$ for the partitioning hypotheses. For all $U \subseteq I$, define
$\boldsymbol{\phi}_U = \tilde{\boldsymbol{\psi}}_U$. This is a valid test of $H_U$ since the partitioning procedure has FWER control. Therefore, we may define a closed testing procedure based on $(\boldsymbol{\phi}_U)_{U \subseteq I}$. Since $\tilde{\boldsymbol{\psi}}_U \leq \tilde{\boldsymbol{\psi}}_V$ whenever $U \subseteq V$, we have
\[
\boldsymbol{\phi}_S^I = \min\{\tilde{\boldsymbol{\psi}}_U\colon S \subseteq U \subseteq I\} = \tilde{\boldsymbol{\psi}}_S,
\]
so the closed testing procedure rejects exactly the same same hypotheses as the partitioning procedure. This proves the second statement.
\end{proof}

\begin{theorem} 
$\mathbf{d}_{\boldsymbol\phi}$ is admissible if and only if the suite $\boldsymbol\phi$ is admissible.
\end{theorem}

\begin{proof} We prove the two counterpoints. Let $\mathbf{d}_{\boldsymbol\phi}$ be inadmissible, and let $\mathbf{\tilde d}$ be a monotone procedure that uniformly improves it. By Theorem \ref{thm sufficiency} there exists $\mathbf{d}_{\boldsymbol\psi} \geq \mathbf{\tilde d}$ that also uniformly improves $\mathbf{d}_{\boldsymbol\phi}$. We have, for every finite $S\subseteq C$,
\[
\boldsymbol\psi_S = \mathds{1}\{\mathbf{d}_{\boldsymbol\psi}^S(S)>0\} \geq \mathds{1}\{\mathbf{d}_{\boldsymbol\phi}^S(S)>0\} = \boldsymbol\phi_S,
\]
Also, by Theorem \ref{thm sufficiency} $\boldsymbol\psi_S$ is a valid local test for $H_S$.

Let $S \subseteq I \subseteq C$, $|I|<\infty$, and $\mathrm{P}\in\Omega$ be such that $\mathrm{P}(E) > 0$ for $E = \{\mathbf{d}_{\boldsymbol\psi}^I(S) > \mathbf{d}_{\boldsymbol\phi}^I(S)\}$. If $E$ happened, by (\ref{def g}) there is a $U \subseteq I$ with $|S\setminus U| = \mathbf{d}_{\boldsymbol\phi}^I(S)$, such that $\boldsymbol\psi_U = 1$ and $\boldsymbol\phi_U=0$. Consequently,
\[
\mathrm{P}(\boldsymbol\phi_U < \boldsymbol\psi_U) \geq \mathrm{P}(E) > 0,
\]
so $\boldsymbol\phi_U$ is inadmissible. Since $|S\setminus U| = \mathbf{d}_{\boldsymbol\phi}^I(S) < \mathbf{d}_{\boldsymbol\psi}^I(S) \leq |S|$, we have $U \neq \emptyset$.

Conversely, let $\boldsymbol\phi_S$ be inadmissible for some finite $\emptyset \subset S \subseteq C$, and let $\boldsymbol\phi_S' \geq \boldsymbol\phi_S$ be a test that uniformly improves it, so that $E = \{\boldsymbol\phi_S' > \boldsymbol\phi_S\}$ has $\mathrm{P}(E) > 0$ for some $\mathrm{P}\in\Omega$.
Define the suite $\boldsymbol\psi$ such that $\boldsymbol\psi_S = \boldsymbol\phi_S'$ and $\boldsymbol\psi_I = \boldsymbol\phi_I$ for $I \neq S$, and consider the monotone procedure $\mathbf{d}_{\boldsymbol\psi}$. Since $\boldsymbol\psi \geq \boldsymbol\phi$ we also have $\mathbf{d}_{\boldsymbol\psi} \geq \mathbf{d}_{\boldsymbol\phi}$. Since $E$ implies $\mathbf{d}_{\boldsymbol\psi}^S(S) > \mathbf{d}_{\boldsymbol\phi}^S(S)$, we see that $\mathbf{d}_{\boldsymbol\psi}$ uniformly improves $\mathbf{d}_{\boldsymbol\phi}$, so the latter is inadmissible.
\end{proof}

\begin{prop} 
If $\mathbf{d}^I \geq \mathbf{d}^I_{\boldsymbol\phi}$ is admissible, then there is an admissible $\boldsymbol\psi$ such that $\mathbf{d}^I = \mathbf{d}^I_{\boldsymbol\psi}$ and, for all $S \subseteq I$, $\boldsymbol\psi_S \geq \boldsymbol\phi^I_{S}$.
\end{prop}

\begin{proof}
Using (\ref{eq trivial scalability}) we may embed $\mathbf{d}^I$ in a monotone procedure. By Theorem \ref{thm sufficiency} we have $\mathbf{d}^I \leq \mathbf{d}^I_{\boldsymbol{\chi}}$, with, for all $S\subseteq I$, $\boldsymbol\chi_S = \mathds{1}\{\mathbf{d}^I(S) > 0\} \geq \mathds{1}\{\mathbf{d}^I_{\boldsymbol{\phi}}(S) > 0\} = \boldsymbol{\phi}^I_S$. For every $S \subseteq I$, if $\boldsymbol{\chi}_S$ is not admissible, let $\boldsymbol\psi_S \geq \boldsymbol{\chi}_S \geq \boldsymbol\phi^I_S$ be a uniform improvement; otherwise, let $\boldsymbol\psi_S = \boldsymbol{\chi}_S$. Without loss of generality we may assume that $\boldsymbol\psi_S$ is admissible. Then $\mathbf{d}^I_{\boldsymbol\psi} \geq \mathbf{d}^I$. Since $\mathbf{d}^I$ is admissible, we have $\mathbf{d}^I_{\boldsymbol\psi} = \mathbf{d}^I$.
\end{proof}

\begin{lemma} 
If $\boldsymbol\phi_S$, $\emptyset\neq S \subseteq I$, is of the form (\ref{eq simeslike}), with $l_{i\mathbin{:}m} \geq l_{i\mathbin{:}n}$ for all $i \geq 1$ and $0 \leq m \leq n$, then
\[
\boldsymbol\phi_S^I = \mathds{1}\{\mathbf{p}_{(i\mathbin{:}S)} \leq l_{i\mathbin{:}\mathbf{h}_I} \textrm{\ for at least one $i=1,\ldots, |S|$} \},
\]
and
\[ 
\mathbf{d}^I(S) = \max_{1\leq u \leq |S|} 1-u+|\{i\in S\colon \mathbf{p}_i \leq l_{u\mathbin{:}\mathbf{h}_I}\}|,
\]
where
\[
\mathbf{h}_I =\max \big\{n\in\{0,\ldots,|I|\}: \mathbf{p}_{(|I|-n+i\mathbin{:}I)}>l_{i\mathbin{:}n}, \textrm{ for } i=1,\ldots,n\big\}.
\]
\end{lemma}

\begin{proof}
Consider first the case $|S| > \mathbf{h}_I$. By definition of $\mathbf{h}_I$, there is an $1 \leq i \leq |S|$ such that
\[
\mathbf{p}_{(i\mathbin{:}S)} \leq \mathbf{p}_{(|I|-|S|+i\mathbin{:}I)}\leq l_{i\mathbin{:}|S|}.
\]
Consequently, $\phi_S =\phi_S^I=1$ for all $S$ with $|S| > \mathbf{h}_I$. Since for such $S$ also $l_{i\mathbin{:}|S|} \leq l_{i\mathbin{:}\mathbf{h}_I}$, the result of the lemma holds if $|S| > \mathbf{h}_I$.

Now consider the case $|S| \leq \mathbf{h}_I$. First, suppose that there exists some $1 \leq i \leq |S|$ with $\mathbf{p}_{(i\mathbin{:}S)} \leq l_{i\mathbin{:}\mathbf{h}_I}$. Take any $V \supseteq S$. If $|V| > \mathbf{h}_I$, we have $\boldsymbol\phi_V=1$ as proved above. If $|V| \leq \mathbf{h}_I$, we have
\[
\mathbf{p}_{(i\mathbin{:}V)} \leq \mathbf{p}_{(i\mathbin{:}S)} \leq l_{i\mathbin{:}\mathbf{h}_I} \leq l_{i\mathbin{:}|V|}
\]
so that $\boldsymbol\phi_V=1$. Since this holds for all $V \supseteq S$, we have $\boldsymbol\phi_S^I=1$.

Next, suppose that there is no $1 \leq i \leq |S|$ with $p_{(i\mathbin{:}S)} \leq l_{i\mathbin{:}\mathbf{h}_I}$. For $j=1, \ldots, |I|$, define $\mathbf{Z}_j$ as a set with $|\mathbf{Z}_j|=j$ such that for all $u \in \mathbf{Z}_j$ and $v \in I \setminus \mathbf{Z}_j$ we have $\mathbf{p}_u \geq \mathbf{p}_v$. Let $\mathbf{W} = S \cup \mathbf{Z}_j$ for some $0\leq j \leq \mathbf{h}_I$ such that $|\mathbf{W}| = \mathbf{h}_I$. If $1 \leq i \leq \mathbf{h}_I-j$, by the assumption we have
\[
\mathbf{p}_{(i\mathbin{:}\mathbf{W})} = p_{(i\mathbin{:}S)} > l_{i\mathbin{:}\mathbf{h}_I} = l_{i\mathbin{:}|\mathbf{W}|}.
\]
If $\mathbf{h}_I-j \leq i \leq \mathbf{h}_I$, since $\mathbf{Z}_j \subseteq \mathbf{Z}_{\mathbf{h}_I}$ we have
\[
p_{(i\mathbin{:}\mathbf{W})} = p_{(i\mathbin{:}\mathbf{Z}_{\mathbf{h}_I})} > l_{i\mathbin{:}\mathbf{h}_I} = l_{i\mathbin{:}|\mathbf{W}|}
\]
because $\boldsymbol\phi_{\mathbf{Z}_{\mathbf{h}_I}}=0$ by definition of $\mathbf{h}_I$. Taken together, this implies that $\boldsymbol\phi_\mathbf{W} = 0$, so $\boldsymbol\phi_\mathbf{W}^I = 0$. This proves the statement about $\boldsymbol\phi_S^I$.

To prove the statement about $\mathbf{d}^I(S)$ we will use
\[
\mathbf{d}^I(S)
= \min_{U \in 2^S} \{|S\setminus U|\colon \boldsymbol\phi_U^I=0\} = \min_{U \in 2^S} \{|I|\colon \boldsymbol\phi_{S\setminus U}^I=0 \}.
\]
As shown above we have $\boldsymbol\phi_{S\setminus U}^I=1$ if and only if for some $1 \leq u \leq |S\setminus U|$ we have $|\{i \in S\setminus U\colon \mathbf{p}_i \leq l_{u\mathbin{:}\mathbf{h}_I}\}| \geq u$, and we may trivially extend the range to $1 \leq u \leq |S|$. Thus, $\boldsymbol\phi_{S\setminus U}^I=0$ if and only if for all such $u$ we have $|\{i \in S\setminus U\colon \mathbf{p}_i \leq l_{u\mathbin{:}\mathbf{h}_I}\}| \leq u-1$. That is, for all such $u$,
\begin{equation} \label{eq_proof_thm1}
|\{i \in S\colon \mathbf{p}_i \leq l_{u\mathbin{:}\mathbf{h}_I}\}| - u + 1 \leq |\{i \in U\colon \mathbf{p}_i \leq l_{u\mathbin{:}\mathbf{h}_I}\}|.
\end{equation}
Denote the left-hand side of (\ref{eq_proof_thm1}) by $\mathbf{g}(u)$ and the right-hand side by $\mathbf{f}(U, u)$. Let $\mathbf{d} = \max_{1 \leq u \leq |S|} \mathbf{g}(u)$. Since $0 \leq \mathbf{d} \leq |S|$ we can pick $\mathbf{U} \subset S$ with $|\mathbf{U}|=d$ such that for all $i \in \mathbf{U}$ and $j \in S\setminus \mathbf{U}$, we have $\mathbf{p}_{i} \leq \mathbf{p}_j$. For this $\mathbf{U}$, for all $1 \leq u \leq |S|$, we have $\mathbf{f}(\mathbf{U},u) \leq |\mathbf{U}| = d$. If $\mathbf{f}(\mathbf{U},u) < d = |\mathbf{U}|$, then $\mathbf{g}(u) \leq \mathbf{f}(S,u) = \mathbf{f}(\mathbf{U},u)$, where the latter step follows by construction of $\mathbf{U}$. If $\mathbf{f}(\mathbf{U},u) = \mathbf{d}$, then $\mathbf{g}(u) \leq \mathbf{d} = \mathbf{f}(\mathbf{U},u)$. We conclude that $\mathbf{U}$ satisfies (\ref{eq_proof_thm1}) for all $1 \leq u \leq |S|$. Obviously, (\ref{eq_proof_thm1}) cannot hold for any $U$ with $|U| < \mathbf{d}$. We conclude that $\mathbf{d}^I(S) = \mathbf{d}$.
\end{proof}

\paragraph*{Proof of equation (\ref{eq HC compare}): $\boldsymbol\phi_S \geq \boldsymbol\psi_S$ for all $S \subseteq I$.}

\begin{proof}
Choose any non-empty $S \subseteq I$. Let $s=|S|$. We have
\begin{eqnarray*}
\boldsymbol\psi_S &=& \mathds{1}\Big\{s-m+\frac{|\{i \in I\colon p_i \leq t\}| - mt - a_m \sqrt{mt(1-t) }}{1-t} > 0 \textrm{\ for some $t \in [0,1)$}\Big\} \\
&=& \mathds{1}\big\{s-m+|\{i \in I\colon p_i \leq t\}| - st - a_m \sqrt{mt(1-t) } > 0 \textrm{\ for some $t \in [0,1)$}\big\} \\
&\leq& \mathds{1}\big\{|\{i \in S\colon p_i \leq t\}| - st - a_s \sqrt{st(1-t) } > 0 \textrm{\ for some $t \in [0,1)$}\big\} \\
&=& \mathds{1}\{\mathbf{f}_S > 0\}\\
& =& \boldsymbol\phi_S
\end{eqnarray*}
\end{proof}

\paragraph*{Proof of equation (\ref{eq interpolation KR}):}
\[
\mathbf{d}^I(S) = 0 \vee \max_{k=1,\ldots,|S|} \big\lceil k - c(1+m\mathbf{p}_{(k\mathbin{:}S)}) \big\rceil,
\]

\begin{proof}
First apply (\ref{def interpolation}) to (\ref{def KR}), yielding
\begin{equation} \label{eq interpolation KR0}
\mathbf{d}^I(S) 
=  0 \vee \max_{i=1,\ldots,m} \big\lceil |\mathbf{K}_i \cap S| - c(1+m\mathbf{p}_{(i\mathbin{:}I)}) \big\rceil.
\end{equation}
To simplify this expression, call $g(i) = |\mathbf{K}_i \cap S| - c(1+m\mathbf{p}_{(i\mathbin{:}I)})$. Let $\boldsymbol\pi_1, \ldots, \boldsymbol\pi_m$ be the permutation such that $\mathbf{K}_i = \{\boldsymbol\pi_1, \ldots, \boldsymbol\pi_i\}$ for all $1\leq i\leq m$. If $\boldsymbol\pi_i \notin S$, then either $i > 1$ and $g(i) \leq g(i-1)$ or $i=1$ and $g(i) < 0$, so we may restrict the maximum in (\ref{eq interpolation KR0}) to values of $i$ with $\boldsymbol\pi_i \in S$. We have $S = \{\boldsymbol\pi_{j_1}, \ldots, \boldsymbol\pi_{j_{|S|}}\}$ with $j_1 < \ldots < j_{|S|}$. If $i=\boldsymbol\pi_{j_k}$, then $g(i) = k - c(1+m\mathbf{p}_{(k\mathbin{:}S)}).$ Therefore, (\ref{eq interpolation KR0}) reduces to (\ref{eq interpolation KR}).
\end{proof}

\begin{lemma} \label{thm adm test} If Assumption \ref{ass shared null} holds, then a statistical test $\boldsymbol\delta$ of hypothesis $H \subseteq \Omega$ is admissible if $\mathrm{P} \in H$ exists such that $\mathrm{P}(\boldsymbol\delta=1) = \alpha$.
\end{lemma}

\begin{proof}
Suppose that that $\mathrm{P}\in H$ exists such that $\mathrm{P}(\boldsymbol\delta=1)=\alpha$, and that $\boldsymbol\delta'$ is a test of $H$ that uniformly improves $\boldsymbol\delta$. We will derive a contradiction under Assumption \ref{ass shared null}. Because $\boldsymbol\delta'$ is a uniform improvement, some $\mathrm{P}' \in\Omega$ exists such that
\begin{equation}
\mathrm{P}'(\boldsymbol\delta'>\boldsymbol\delta)>0. \label{eq theta prime}
\end{equation} By Assumption \ref{ass shared null}, (\ref{eq theta prime}) remains valid if we replace $\mathrm{P}'$ by $\mathrm{P}$. Consequently, since $\boldsymbol\delta'\geq\boldsymbol\delta$, and since $\{\boldsymbol\delta'>\boldsymbol\delta\}$ and $\{\boldsymbol\delta=1\}$ are disjoint, we have
\[
\mathrm{P}(\boldsymbol\delta'=1) = \mathrm{P}(\boldsymbol\delta'>\boldsymbol\delta) + \mathrm{P}(\boldsymbol\delta=1) > \alpha
\]
which contradicts that $\boldsymbol\delta'$ is a valid test of $H$.
\end{proof}

\begin{lemma} 
$\mathbf{d}^I$ is consonant and coherent if and only if, for every disjoint $V,W \subseteq I$,
\begin{equation} \label{eq d consonant}
\mathbf{d}^I(V \cup W) = \mathbf{d}^I(V) + \mathbf{d}^I(W).
\end{equation}
\end{lemma}

\begin{proof}
Choose any $V,W \subseteq I$ disjoint. Call $S = V \cup W$. We use complete induction on $|S|$. Suppose that (\ref{eq d consonant}) holds for all sets smaller than $S$. If $V=\emptyset$ or $W=\emptyset$ the result is trivial, so we assume and $V,W \neq \emptyset$. If $\mathbf{d}^I(S) = 0$ the result follows immediately from Lemma \ref{thm coherence}, so we may assume $\mathbf{d}^I(S) > 0$. By consonance there is an $i \in S$ such that $\mathbf{d}^I(\{i\}) = 1$. Without loss of generality, suppose that $i \in W$. By Lemma \ref{thm coherence} and the induction hypothesis we have
\[
\mathbf{d}^I(S)\ =\ \mathbf{d}^I(V \cup W \setminus \{i\}) + \mathbf{d}^I(\{i\})\ =\
\mathbf{d}^I(V) + \mathbf{d}^I(W \setminus \{i\}) + \mathbf{d}^I(\{i\}).
\]
Since $V\neq\emptyset$, $W \subset S$ and we may use the induction hypothesis once more, saying that $\mathbf{d}^I(W \setminus \{i\}) + \mathbf{d}^I(\{i\}) = \mathbf{d}^I(W)$ to obtain (\ref{eq d consonant}).

For the converse, suppose that (\ref{eq d consonant}) holds. By Lemma \ref{thm coherence}, $\mathbf{d}^I$ is coherent. Choose some $\emptyset\neq S\subseteq I$ such that $\mathbf{d}^I(S)>0$. We use complete induction on $|S|$ to show that $\mathbf{d}^I(\{i\})=1$ for some $i \in S$. Suppose the result holds for all sets smaller than $S$. Choose any $i \in S$ and let $V=\{i\}$ and $W=S\setminus\{i\}$. By (\ref{eq d consonant}), either $\mathbf{d}^I(V)>0$ or $\mathbf{d}^I(W)>0$. If the former, we have the result immediately. If the latter, we have use the induction hypothesis to conclude that $\mathbf{d}^I(\{i\})=1$ for some $j \in W \subset S$.
\end{proof}

\begin{lemma} 
If $\mathbf{d}^I_{\boldsymbol\phi}$ is consonant, $\mathbf{d}^I_{\boldsymbol\phi} = \mathbf{r}^I_{\boldsymbol\phi}$; otherwise $\mathbf{d}^I_{\boldsymbol\phi}$ uniformly improves $\mathbf{r}^I_{\boldsymbol\phi}$.
\end{lemma}

\begin{proof}
Choose $S \subseteq I$ en let $\mathbf{V} = S \cap \mathbf{R}^I_{\boldsymbol\phi}$ and $\mathbf{W} = S \setminus \mathbf{R}^I_{\boldsymbol\phi}$. By definition of $\mathbf{R}^I_{\boldsymbol\phi}$, for every $i \in \mathbf{V}$ we have $\boldsymbol\phi^I_{\{i\}} = 1$, so $\boldsymbol\phi^I_{J} = 1$ for all $\emptyset \neq J \subseteq \mathbf{V}$, so $\mathbf{d}^I_{\boldsymbol\phi}(\mathbf{V}) = |\mathbf{V}|$ by (\ref{def d}).

If $\mathbf{d}^I_{\boldsymbol\phi}$ is consonant,
by definition of $\mathbf{R}^I_{\boldsymbol\phi}$, for every $i \in \mathbf{W}$ we have $\boldsymbol\phi^I_{\{i\}} = 0$. By consonance, we must have $\boldsymbol\phi^I_{\mathbf{W}} = 0$, so $\mathbf{d}^I_{\boldsymbol\phi}(\mathbf{W}) = 0$  by (\ref{def d}). By Lemma \ref{thm VW consonant} we have $\mathbf{d}^I_{\boldsymbol\phi}(S) = |\mathbf{V}|$, which proves the first statement.

If $\mathbf{d}^I_{\boldsymbol\phi}$ is not consonant, by Lemma \ref{thm coherence} we have $\mathbf{d}^I_{\boldsymbol\phi}(S) \geq \mathbf{d}^I_{\boldsymbol\phi}(\mathbf{V}) = \mathbf{r}^I_{\boldsymbol\phi}(S)$. Moreover, $S\subseteq I$ exists such that for some $\mathrm{P}\in\Omega$ we have with positive probability that $\mathbf{d}^I_{\boldsymbol\phi}(S)>0$ and $S \cap \mathbf{R}^I_{\boldsymbol\phi} =\emptyset$, or $\mathbf{d}^I_{\boldsymbol\phi}$ would be consonant. For this $S$ with positive probability $\mathbf{d}^I_{\boldsymbol\phi}(S) > 0 = \mathbf{r}^I_{\boldsymbol\phi}(S)$.
\end{proof}

\begin{prop} 
If $\mathbf{R}$ is admissible, then a consonant suite $\boldsymbol\phi$ exists such that $\mathbf{R} = \mathbf{R}_{\boldsymbol\phi}$.
\end{prop}

\begin{proof}
For all finite $S \subseteq C$ consider $\boldsymbol\phi_S = \mathds{1}\{\mathbf{R}^S \neq \emptyset\}$. Then for all finite $S \subseteq I \subseteq C$
\begin{equation} \label{eq local test FWER}
\boldsymbol\phi_S^I = \mathds{1}\{\textrm{$\mathbf{R}^J \neq \emptyset$ for all $S \subseteq J \subseteq I$}\}.
\end{equation}

We show that $\boldsymbol\phi$ is consonant, i.e.\ if $\boldsymbol\phi^I_S = 1$, there is an $i \in \mathbf{R}_S \neq \emptyset$ such that $\boldsymbol\phi^I_{\{i\}} = 1$. We proceed by complete induction on $|S|$ downward from $|I|$. Assume that for all $V \subseteq I$ with $|V| > |S|$, it holds that $\boldsymbol\phi^I_V = 1$ implies that an $i \in \mathbf{R}_V\neq \emptyset$ exists such that $\boldsymbol\phi^I_{\{i\}} = 1$, but the same does not hold for $V=S$. We will derive a contradiction. Since $\boldsymbol\phi^I_S = 1$, indeed, by (\ref{eq local test FWER}), $\mathbf{R}^S \neq \emptyset$. Choose $i \in \mathbf{R}^S$ and $J \ni i$. If $J \subseteq S$, we have $i \in \mathbf{R}_S \cap J \subseteq \mathbf{R}^J$ by monotonicity, so $\mathbf{R}^J \neq \emptyset$. If not $J \subseteq S$, we have $\boldsymbol\phi^I_{J \cup S} = 1$, since $\boldsymbol\phi^I_{S} = 1$ and $S \subset J \cup S$. By the induction hypothesis there exists a $j \in \mathbf{R}^{J \cup S}$ such that $\boldsymbol\phi^I_{\{j\}} = 1$. Now either $j \in S$ or $j \in J \setminus S$. In the former case $j \in S \cap \mathbf{R}^{J \cup S} \subseteq \mathbf{R}^S$ by monotonicity, and we have a contradiction. Therefore, we must have $j \in J \setminus S$, in which case $j \in \mathbf{R}^{J \cup S} \cap J \subseteq \mathbf{R}^J$, also by monotonicity, so $\mathbf{R}^J \neq \emptyset$. Since $\mathbf{R}^J \neq \emptyset$ for all $J \ni i$, we have $\boldsymbol\phi^I_{\{i\}} = 1$ by (\ref{eq local test FWER}). This proves consonance.

Since $\mathbf{R}$ is monotone, we have
\[
\boldsymbol\phi_S^I \geq \mathds{1}\{\textrm{$\mathbf{R}^I \cap J \neq \emptyset$ for all $S \subseteq J \subseteq I$}\} = \mathds{1}\{\mathbf{R}^I \cap S \neq \emptyset\}.
\]
Clearly, by (\ref{eq closed FWER}), $\mathbf{R}^I_{\boldsymbol\phi} \geq \{i \in I\colon \mathbf{R}^I \cap \{i\} \neq \emptyset\} = \mathbf{R}^I$ for all finite $I \subseteq C$. Since $\mathbf{R}$ is admissible we must have $\mathbf{R}^I_{\boldsymbol\phi} = \mathbf{R}^I$.
\end{proof}

\begin{prop} 
If $\mathbf{R}$ is admissible, then an admissible suite $\boldsymbol\psi$ exists such that $\mathbf{R} = \mathbf{R}_{\boldsymbol\psi}$.
\end{prop}

\begin{proof}
Let $\mathbf{R}$ be admissible. By Proposition \ref{thm Romano}, $\mathbf{R} = \mathbf{R}_{\boldsymbol\phi}$, with $\boldsymbol\phi$ consonant. If $\boldsymbol\phi$ is admissible, we are done. Otherwise, let $\boldsymbol\psi$ uniformly improve $\boldsymbol\phi$. Without loss of generality, we may assume that $\boldsymbol\psi$ is admissible. Then $\mathbf{R}_{\boldsymbol\psi} \geq \mathbf{R}_{\boldsymbol\phi}$, and we must have equality since $\mathbf{R}$ is admissible.
\end{proof}

\begin{prop} 
If $\boldsymbol\phi$ is consonant and admissible, then $\mathbf{R}_{\boldsymbol\phi}$ is admissible.
\end{prop}

\begin{proof}
Let $\boldsymbol\phi$ be consonant. By Lemma \ref{thm r inadmissible}, we have $\mathbf{d}_{\boldsymbol\phi} = \mathbf{r}_{\boldsymbol\phi}$. Suppose $\mathbf{R}$ uniformly improves $\mathbf{R}_{\boldsymbol\phi}$. Then $\mathbf{r} = (\mathbf{r}^I)_{I \subseteq C, |I|<\infty}$ with $\mathbf{r}^I(S) = |S \cap \mathbf{R}^I|$ uniformly improves $\mathbf{r}_{\boldsymbol\phi} = \mathbf{d}_{\boldsymbol\phi}$, so $\boldsymbol\phi$ is inadmissible by Theorem~\ref{thm admissibility}.
\end{proof}

\section{A sufficient condition for using bisection to calculate $\mathbf{h}_I$}

In this section we show that $\mathbf{h}_I$ can be calculated by bisection in $O(m\log(m))$ time if, for all $2\leq i < n$, we have
\begin{equation} \label{eq condition bisection}
l_{i-1\mathbin{:}n-1} \leq l_{i\mathbin{:}n}
\end{equation}

By Lemma \ref{thm shortcut}, we have
\[
\mathbf{h}_I =\max \big\{s\in\{0,\ldots,|I|\}: \mathbf{p}_{(|I|-s+i)}>l_{i\mathbin{:}s}, \textrm{ for } i=1,\ldots,s\big\}.
\]
For bisection to be used we need that the condition that
\begin{equation} \label{eq bisection induction}
\mathbf{p}_{(|I|-s+i)}>l_{i\mathbin{:}s}, \textrm{ for all } i=1,\ldots,s
\end{equation}
holds for all $s \leq \mathbf{h}_I$. This follows immediately by induction, since (\ref{eq bisection induction}) holds for $s = \mathbf{h}_I$, and if (\ref{eq bisection induction}) holds, then for any $1 \leq i \leq s-1$, we have
\[
\mathbf{p}_{(|I|-(s-1)+i)} = \mathbf{p}_{(|I|-s+(i+1))} > l_{i+1\mathbin{:}s} \geq l_{i\mathbin{:}s-1},
\]
so (\ref{eq bisection induction}) holds for the next lower value of $s$.

It is easy but tedious to check that (\ref{eq condition bisection}) holds for (\ref{eq l HC}), but it does not hold for (\ref{eq KR local}) or (\ref{eq KR improved critical}).

\section*{Acknowledgements}

This paper was inspired by many discussions during the workshop ``Post-selection Inference and Multiple Testing'' in Toulouse, Februari 2018, organized by Pierre Neuvial, Etienne Roquain and Gilles Blanchard. We thank the organizers and all participants, and especially Ruth Heller for asking the question that triggered this research project. We thank Jonathan Rosenblatt and the Israel Science Foundation for financing the computing equipment used for the simulations (grants 926/14 and 900/16). Jelle Goeman was supported by NWO VIDI grant 639.072.412.

\bibliographystyle{apalike}
\bibliography{sufficiency}

\end{document}